%% file: main.tex
\documentclass{article}
\usepackage{amsmath}
\usepackage{amssymb}
\usepackage[margin=1in]{geometry}
\usepackage[utf8]{inputenc}
\usepackage{graphicx}
\usepackage{caption}
\usepackage{subcaption}
\usepackage{float}
\usepackage{setspace}
\usepackage[titletoc,title]{appendix}
\usepackage{xcolor}
\allowdisplaybreaks
\usepackage[noblocks]{authblk}
\usepackage{amsthm}
\usepackage{tikz}
\usetikzlibrary{shapes,arrows,chains}
\usepackage{mathtools}
\usepackage{pgfplots}
\usepackage{soul}
\usepackage[style=ext-authoryear-comp]{biblatex}
\usepackage{multirow}

\newcommand{\bb}{b}
\newcommand{\A}{A}
\newcommand{\erlangshape}{n}
\newcommand{\reservationprofit}{Z}

\newtheorem{proposition}{Proposition}

\newtheorem{lemma}{Lemma}
\newtheorem{corollary}{Corollary}

\newlength\fwidth

\usepackage{xpatch}
\xpatchcmd{\author}{\relax#1\relax}{\relax\detokenize{#1}\relax}{}{}

\title{Contingent Penalty and Contingent Renewal Supply Contracts in High-Tech Industry}
\author{Mirjam S. Meijer \thanks{Corresponding author. Email: m.s.meijer@tue.nl}}
\author{Willem van Jaarsveld}
\author{Ton de Kok}
\affil{Department of Industrial Engineering and Innovation Sciences, Eindhoven University of Technology, Eindhoven, The Netherlands}
\author{Christopher S. Tang}
\affil{UCLA Anderson School of Management, University of California, Los Angeles, Los Angeles, California, USA}
\date{}

\newcommand{\expect}{\operatorname{\mathbb{E}}\expectarg}

\DeclarePairedDelimiterX{\probarg}[1]{(}{)}{%
  \ifnum\currentgrouptype=16 \else\begingroup\fi
  \activatebar#1
  \ifnum\currentgrouptype=16 \else\endgroup\fi
}
\DeclarePairedDelimiterX{\expectarg}[1]{[}{]}{%
  \ifnum\currentgrouptype=16 \else\begingroup\fi
  \activatebar#1
  \ifnum\currentgrouptype=16 \else\endgroup\fi
}
\newcommand{\activatebar}{%
  \begingroup\lccode`\~=`\|
  \lowercase{\endgroup\let~}\innermid
  \mathcode`|=\string"8000
}

\newcommand{\comment}[1]{}

\begin{document}
\maketitle
\begin{abstract}
\noindent
{\bf{Abstract:}} Unlike consumer goods industry, a high-tech manufacturer (OEM) often amortizes new product development costs over multiple generations, where demand for each generation is based on advance orders and additional uncertain demand. Also, due to economic reasons and regulations, high-tech OEMs usually source from a single supplier. Relative to the high retail price, the wholesale price for a supplier to produce high-tech components is low. Consequently, incentives are misaligned: the OEM faces relatively high under-stock costs and the supplier faces high over-stock costs. 

In this paper, we examine supply contracts that are intended to align the incentives between a high-tech OEM and a supplier so that the supplier will invest adequate and yet non-verifiable capacity to meet the OEM's uncertain demand. When focusing on a single generation, the manufacturer can coordinate a decentralized supply chain and extract all surplus by augmenting a traditional wholesale price contract with a ``contingent penalty'' should the supplier fail to fulfill the OEM's demand.  When the resulting penalty is too high to be enforceable, we consider a new class of ``contingent renewal'' wholesale price contracts with a stipulation: the OEM will renew the contract with the incumbent supplier for the next generation {\em{only when}} the supplier can fulfill the demand for the current generation.  By using non-renewal as an implicit penalty, we show that the contingent renewal contract can coordinate the supply chain. While the OEM can capture the bulk of the supply chain profit, this innovative contract cannot enable the OEM to extract the entire surplus.
\end{abstract}

\section{Introduction}\label{sec:introduction}
This paper examines a supply contract problem arising from a high-tech supply chain that involves an original equipment manufacturer (OEM) who designs and manufactures state-of-the-art systems and a focal supplier who makes a critical component for the system. 
Examples of OEMs in the high-tech industry include commercial aircraft (Boeing, Airbus), medical imaging equipment (Philips, Siemens, GE), and lithography systems for the semiconductor industry (ASML, CANON). At the same time, suppliers who make critical components include Yuasa (who makes the batteries for Boeing's 787), VDL/ETG (who makes the wafer handler for ASML), and Neways (who makes control systems for Philips). To reduce development cost and development time in the high-tech industry, most OEMs act as system ``integrators'': they initiate and develop different generations of a product, which they engineer together with hundreds of suppliers \parencite{tang2009managing}.
However, based on our discussion with ASML in the Netherlands, we learned that, unlike mass produced consumer products such as apparel and home furnishings, the sourcing and supplier management for high-tech components are fundamentally different as follows:

\begin{enumerate}
    \item {\bf{Multiple product generations}}.  Because the research and development cost of a new product involves billions of euros, each new product has multiple generations so that the OEM can continue to improve the design for each generation instead of starting from scratch.  Therefore, accounting for the life of different generations, the life cycle of a new product can last for decades.  For example, Boeing's 747 was launched in 1970 and retired in 2018.  
    \item {\bf{Advance orders and a highly uncertain additional demand}}.  Like any high-tech products, the demand for a new product (or a new generation of the product) is highly uncertain because some customers are risk-averse in adopting new technologies especially when they are uncertain about product performance.  To reduce demand uncertainty, many OEMs (e.g., Boeing, ASML) encourage customers to place their orders in advance even though some customers prefer to order after product launch. For example, Boeing received some 900 advance orders for the Boeing 787, but many airlines only ordered after the aircraft was in production \parencite{nolan09}.  Hence, the \emph{base demand} associated with the advance orders is known, but a highly uncertain \emph{additional demand} remains prevalent. 
    \item {\bf{Single-sourced.}} Because the production of high-tech components requires component-specific equipment and technology-specific technical staff, the production \emph{capacity} of a high-tech component is also generation-specific.\footnote{This is because new technology is involved for each new generation, the extended use of equipment for older generation production is normally infeasible or impractical.}  Also, due to high demand uncertainty, suppliers are reluctant to participate unless it is sole-sourced.  For these reasons, the OEM would normally work with the same supplier for designing and producing components for multiple product generations.\footnote{Also, due to the underlying advanced technology, OEMs usually work with a single supplier for the development of each component to foster close cooperation with the intention of a longer-term relationship. Examples include Boeing with Alcoa, ASML with VDL/ETG, Philips with Neways, etc.}  
    \item {\bf{Non-verifiable supply capacity.}}  While the OEM can audit the equipment acquired by the supplier for a specific generation, the actual production capacity is difficult for the OEM to verify. Production often involves technical staff and engineers to operate the required equipment, but the available staff can also be assigned to work on other OEMs' orders. Therefore, the actual production capacity is difficult to judge for the OEM. While the actual capacity is not verifiable ex-ante, the OEM will know the supplier's capacity only when the supplier failed to meet the ex-post realized demand (placed by the OEM).
    \item  {\bf{Very high under-stock cost.}}  In the high-tech industry, the selling price of each system is in millions of euros. At the same time, the research and development cost is in billions of euros, while the unit cost of a component is relatively very low: it can range from a few hundred euros to tens of thousands of euros. 
\end{enumerate}

Other issues such as varying product quality \parencite[see e.g.][]{transchel16} may also play a role, but are outside the scope of this paper.
The above context creates the following challenges faced by many OEMs (e.g., ASML) in the high-tech industry.  First, due to high selling price and low unit production cost (excluding the research and development cost), the OEM incurs a very high under-stock cost and would like the supplier to invest ample capacity to meet both the base demand and additional demand.   However, because capacity is very costly and generation-specific and the wholesale price is relatively low, the supplier incurs high over-stock cost and has little incentive to invest in ample capacity.  Because of the underlying misaligned incentives, many high-tech OEMs face shortages of components to meet uncertain demand.\footnote{For example, a shortage of aircraft-grade fasteners, allegedly resulting from Alcoa's insufficient capacity investments, caused headline delivery delays of the Boeing 787 \parencite{alcoasupplydisaster}.}  Second, due to non-verifiable supply capacity ex-ante, the OEM cannot contract directly based on the reserved capacity as examined in the literature \parencite{brown2003impact,roels2017win}. 

These two challenges motivate us to develop a new class of supply contracts that is intended to entice the supplier to invest sufficient capacity to coordinate the decentralized supply chain by aligning the incentives of the OEM and the supplier when the capacity is not verifiable ex-ante.  In addition to the known base demand, we assume that the additional uncertain demand follows an exponential distribution to ensure  tractability. We first examine the classic wholesale price contract arising from a single product generation that may occur when the OEM treats the sourcing decision of different product generations separately.  When the OEM pays the supplier a wholesale price for each unit delivered, we re-confirm a well-known result: a traditional wholesale price contract can coordinate the supply chain only when the wholesale price equals the OEM's profit margin, which the OEM will not oblige \parencite{cach03}.  However, when the OEM offers a wholesale price while imposing a ``shortfall penalty'' (i.e., the supplier fails to meet the OEM's demand), we find that the ``augmented'' wholesale price contract can enable the OEM to coordinate the supply chain while capturing all the profit.  Hence, such an augmented contract is optimal to the OEM for managing the supply contract for a single product generation.  

While the augmented wholesale price contract is optimal and can coordinate the decentralized supply chain for a single product generation, it may not be practical when the penalty is too high to be enforced upon financially constrained suppliers. Our analysis reveals that this happens primarily when the cost of underage is extremely high, a situation that occurs in the high-tech setting as selling prices are much higher than the supplier's capacity and production costs. 

In those cases, instead of imposing a direct shortfall penalty, we consider a new class of supply contracts that takes the sourcing of multiple product generations into consideration. 
Specifically, we consider a class of ``contingent renewal'' wholesale price contracts that can be described as follows: in addition to the wholesale price, the OEM will renew the contract with the incumbent supplier for the production of the next generation only when it has sufficient capacity to fulfill the OEM's demand for the current generation.  Observe that this contingent contract creates an ``indirect penalty'' associated with non-renewal that affects the supplier's future profit,\footnote{Discussions with high-tech manufacturers have revealed that long-term cooperations (and the possibility to terminate them) are key to entice suppliers to comply.} and this contingent contract is enforceable because the OEM has the option to work with different suppliers for different generations (especially when the supplier capacity is generation-specific). 
By expanding our model to capture the characteristics of multiple product generations, we find that the contingent renewal supply contract can enable the OEM to coordinate the decentralized supply chain; however, the OEM cannot extract the entire surplus from the supplier. Even so, we find that the OEM can capture the bulk of the total profit of the entire supply chain when the selling price is much higher than the supplier's capacity and production costs and when the supplier has a high valuation of future profits or when a substantial fraction of total demand is ordered in advance. 

This paper is organized as follows.  We review relevant literature in \S\ref{sec: literature review}.  In \S\ref{sec:singleproductgeneration}, we focus our analysis for the single product generation case and show that, by augmenting the wholesale price contract with a contingent penalty, the OEM can coordinate the supply chain optimally by extracting the entire surplus from the supplier.  \S\ref{sec: multi generation} extends our analysis to the multi-generation case.  By developing a wholesale price contract with endogenous renewal probability, the OEM can coordinate the supply chain but it cannot extract the entire surplus from the supplier.  However, the coordinating contingent wholesale price contract enables the OEM to capture the bulk of the total supply chain profit under certain conditions. In \S\ref{sec: extensions} we show that our results continue to hold when we consider for example other demand distributions. The paper concludes in \S\ref{sec: conclusions}. All proofs are provided in the Appendix. 

\section{Literature Review}\label{sec: literature review}
Some of the contracts studied in this paper may be renewed over multiple periods.  
\textcite{tay07,taylor2007simple} study informal, relational contracts that create incentives for the supplier to invest in production capacity repeatedly over multiple periods.
\textcite{tay07} derive the optimal self-enforcing relational contract. However, this contract may be complex and therefore difficult to implement, for which reason an intermediate-complexity relational contract that performs well in many parameter settings is considered as well.
Since these contracts are often still difficult to implement in practice, \textcite{taylor2007simple} consider simple relational contracts that consist of agreements on price only or on price and quantity. They specify conditions under which a price-only contract is most applicable and when a price-quantity contract is more efficient and compare the performance of the most suitable one to that of the optimal relational contract from \textcite{tay07}. They show that for large discount factors the performance of the simple relational contracts is close to that of the optimal relational contract, but for moderate capacity cost and discount factors the loss compared to the optimal relational contract is substantial.
Similar to \textcite{tay07,taylor2007simple}, \textcite{sun2014} show that informal long-term relationships can be sustained when the discount factor of future profits is sufficiently high, taking into account the effect of turbulent markets.

While Taylor and Plambeck focus on comparing various relational contracts, we focus on assessing the relative merits of simple contracts and their applicability inspired by a practical application in high-tech supply chains.  We show that wholesale price contracts with penalty can coordinate the supply chain but may not always be enforceable.  Although our renewable wholesale price contracts are similar to the price-only relational contracts, our contracts perform better especially when the valuation of future profits is high (cf. Taylor and Plambeck). 

Our model explicitly distinguishes between a known base demand (known before product launch) and an uncertain tail demand (revealed after product launch). We derive new analytic results for this model to generate clear insights that are particularly relevant in the high-tech setting. In this setting retail prices are typically high, and for the corresponding asymptotic limit we derive a simple closed-form expression of the division of profit over the OEM and supplier. In this limit, the renewable wholesale price contracts result in the OEM capturing a large share of the profit, which motivates OEMs to adopt this very simple contract in practice. Additionally, we show that a higher known base demand results in a higher share of profit captured by the manufacturer.  

Our paper is also related to the capacity reservation literature \parencite[see e.g.][]{barnes02,brown2003impact,erkoc05,ren10,roels2017win}.  Unlike these single-period models, our high-tech industry context lends itself to long-term partnerships involving repeated interactions with the supplier \parencite{jin07}.
While \textcite{serel07} considers a multi-period capacity reservation contract between a manufacturer and a long-term supplier, we consider the case when capacity cannot be verified and when the contract renewal hinges on the supplier's performance.  
Long-term supply contracts are also considered by \textcite{frascatore08}. \textcite{frascatore08} conclude that the supplier can be induced to create higher capacity by offering a contract that spans multiple periods and that including a penalty can encourage the supplier to invest in the supply chain optimal capacity. 
Our model differs significantly from the setting considered by \textcite{frascatore08}, as they consider a relation spanning a fixed number of periods while we consider a potentially infinite collaboration where continuation depends on the supplier's investment decisions. Furthermore, we show that direct penalties may not be enforceable in the considered setting.
Decisions on continuing a supply chain relationship are studied by \textcite{pfeiffer10} in a setting with information asymmetry. It is found that the threat of switching between suppliers can be used as an instrument to reduce information asymmetry.

Our model is related to Vendor Managed Inventory (VMI) programs in which the supplier is responsible for replenishing the manufacturer's inventory. The past decades VMI has been investigated increasingly \parencite[see e.g.][]{fry01,lee18,hu18}.  
\textcite{corb01} and \textcite{chintapalli17} study how supply chain efficiency can be increased by delegating replenishment decisions to the supplier, as inefficiency due to sub-optimal order quantities is reduced.
However, these studies focus on one-off interactions or are based on deterministic demand. 
\textcite{guan10} consider repeated interaction in a stochastic demand setting and find that parties are  more willing to share private information when engaging in repeated interactions.
Our study is fundamentally different from this stream of work. In our context, the contract renewal probability is endogenously dependent on the supplier's capacity decision and the contract will not be renewed should the supplier fail to fulfill the OEM's uncertain demand.

Finally, the study by \textcite{sieke12} is closely related to our study. They study supply chain coordination using service level contracts that enforce pre-specified service levels with financial penalty payments. It is concluded that for the considered types of contracts for every service level there exists a contract that coordinates the supply chain, assuming there are no limitations on the penalty costs. We show that in practice there often are limitations on the penalty costs and suggest an alternative type of contracts that calls for non-renewal should the supplier fail to fulfill the OEM's uncertain demand.

\section{Contracting for one product generation}\label{sec:singleproductgeneration}
In this section, we consider a contracting issue arising from the development of one product generation between the OEM and the supplier.  This setting occurs when the underlying product has only one generation or when the OEM deals with the supply contract for different product generations separately (i.e., one generation at a time). However, in \S4 we shall extend our analysis to the case when the OEM takes the development of multiple product generations into consideration when contracting with a supplier.

\subsection{Centralized Supply Chain: a benchmark}
To establish a benchmark, let us begin by analyzing a supply chain under centralized control in which the OEM and the supplier are managed by a central body.  The centralized supply chain makes the  
capacity decision $x$ ``before'' demand $D$ is realized. The cost for investing in the requisite capacity is $c\geq 0$ per unit, and the demand (for the single production generation) is $D=\bb+\A$, where $\bb\geq 0$ is the ``base demand" that is known in advance (e.g., advance orders) and $\A$ represents the uncertain ``additional'' demand that is realized after the product is launched, with $A\sim Exp(\lambda)$.\footnote{The decomposition of demand into advance orders and uncertain late orders is for example seen in the case of Boeing, where over 900 advance orders were received before the Boeing 787 was actually produced, while some airlines would place the order only after the new aircraft is in production \parencite{nolan09}.} 
To capture the notion as explained in \S1 that many customers are reluctant to place their orders of high-tech products in advance, we shall assume that the expected additional demand is greater than the base demand so that $\expect{\A} = \frac{1}{\lambda} \ge \bb$. 

Given the capacity $x$ established at $c$ per unit  ``before'' the demand is realized, the centralized supply chain can produce up to $x$ units at cost $k$ per unit ``after'' the demand $D$ is realized. 
If $D \le x$, then the gross revenue is equal to $(r-k)D$, where $r$ is the exogenously given retail price. To capture the notion that high tech products have a high net profit margin $(r-k-c)$, we shall assume that  $r-k-c > c$.\footnote{Without the assumption that $r-k-c>c$, additional conditions are necessary in Proposition \ref{prop:ws} to satisfy the supplier's participation constraint. All other results presented in this paper continue to hold without this assumption.\label{fn:assumption r-k-c}}

By assuming that unmet demand is lost and the established capacity $x$ has no salvage value,{\footnote{Our model can be extended to the case when capacity has salvage value.  For ease of exposition, we scale the salvage value to zero.}} the profit of the centralized supply chain $\Pi(x)$ can be written as: 
\begin{equation}\label{eq:sp_sc_p}
    \Pi(x) = -cx + \expect{(r-k)\min\{D, x\}} =  -cx + \expect{(r-k)\min\{\bb+\A,x\}}
\end{equation}

By considering the first order condition and by using the fact that $\A\sim Exp(\lambda)$ along with our assumptions
that $r - k - c > c$ and $\frac{1}{\lambda} \ge \bb$, we get: 
\begin{proposition}\label{prop:sp_centralized}
When the supply chain is controlled centrally, the optimal capacity $x^* = \bb+ \frac{1}{\lambda}\ln\left(\frac{r-k}{c}\right)$ and the corresponding optimal supply chain profit $\Pi^*=(r-k-c)\left(\bb+\frac{1}{\lambda}\right) - \frac{c}{\lambda}\ln\left(\frac{r-k}{c}\right) > 0$. 
\end{proposition}

\noindent
Besides the base demand $\bb$, the optimal capacity $x^*$ includes some ``additional capacity'' $\frac{1}{\lambda}\ln\left(\frac{r-k}{c}\right) > 0$  that is intended to satisfy the uncertain additional demand $\A$.  This additional capacity is increasing in $\expect{\A} = \frac{1}{\lambda}$ and the gross margin $(r-k)$, but it is decreasing in the capacity unit cost $c$. Armed with the ``first best solution'' $x^*$ as a benchmark, we now consider the case when the underlying supply chain is decentralized. 

\subsection{Decentralized Supply Chain: wholesale price contracts}\label{sec:wholesale}

In a decentralized supply chain, the OEM delegates the capacity investment decision to an external supplier and offers a wholesale contract to the supplier based on a wholesale price $w$ (a decision variable). After the OEM has decided on and communicated $w$, the supplier decides on capacity $x$, where $x$ is ``not verifiable'' by the OEM as explained in \S\ref{sec:introduction}.  In this case, the supplier establishes capacity $x$ at cost $c$ per unit before the demand $D$ is realized.  After the demand is realized, the supplier can produce up to $x$ at cost $k$ per unit.  Hence, the supplier's profit for any given wholesale price $w$ is: 
\begin{equation}\label{eq:dc_ps}
    \tilde{\pi}_s(x,w) = -cx + \expect{(w-k)\min\{D, x\}} =  -cx + \expect{(w-k)\min\{\bb+\A,x\}}
\end{equation}
By considering the first order condition and the fact that $\A\sim Exp(\lambda)$, the supplier's optimal capacity is:
\begin{equation}\label{eq:sp_dc_x}
\tilde{x}(w)=\bb+\frac{1}{\lambda} \ln\left(\frac{w-k}{c}\right)
\end{equation}
if and only if $w \ge c+k$.\footnote{Clearly, if the wholesale price is below cost (i.e., $w<c+k$), $\tilde{x}(w)=0$.} Substituting \eqref{eq:sp_dc_x} into \eqref{eq:dc_ps}, the corresponding optimal supplier profit is:
\begin{equation} \label{pisd}
    \tilde{\pi}_s(w) \equiv \tilde{\pi}_s(\tilde{x}(w),w)=(w-c-k)\left(\bb+\frac{1}{\lambda}\right)-\frac{c}{\lambda} \ln\left(\frac{w-k}{c}\right).
\end{equation}
The following lemma asserts that the supplier's participation constraint (i.e., $\tilde{\pi}_s(w) \ge 0$) holds\footnote{To ease our exposition and without loss of generality, we scale the minimum acceptable profit for the supplier to accept a contract to $0$, which is standard approach being used in the supply contract literature (e.g., \cite{cach03}).  For completeness, we shall extend our analysis to the case when this minimum acceptable profit is any $Z \ge 0$ in \S 5.1.} when the wholesale price is greater than the effective unit cost $c + k$. Hence, the supplier participation constraint holds when $w \ge c+k$.
\begin{lemma}
The supplier's profit $\tilde{\pi}_s(w) = (w-c-k)\left(\bb+\frac{1}{\lambda}\right)-\frac{c}{\lambda} \ln\left(\frac{w-k}{c}\right) \ge 0$ if $w \ge k + c$. \label{lem:helperlemma} 
\end{lemma}

By comparing $\tilde{x}(w)$ given in (\ref{eq:sp_dc_x}) and $x^*$ given in Proposition 1, we get: 
\begin{corollary}\label{cor:ws_c}
Under the wholesale price contract, one can coordinate the decentralized supply chain so that $\tilde{x}(w) = x^*$ if and only if $w=r$. 
\end{corollary}
\noindent
In view of Corollary \ref{cor:ws_c}, the OEM can coordinate the supply chain by setting $w = r$ to entice the supplier to reserve the capacity $\tilde{x}(r) = x^*$; however, the OEM will earn nothing which is undesirable.  

\subsubsection{Optimal wholesale price}  

Instead of focusing on supply chain coordination by setting $w = r$, let us consider the case when the OEM optimizes its own profit. 

By anticipating the supplier's optimal capacity investment $\tilde{x}(w) = \bb+\frac{1}{\lambda} \ln\left(\frac{w-k}{c}\right)$ as given in (\ref{eq:sp_dc_x}),
the OEM's profit associated with any $w \ge c + k$ (to ensure supplier participation) is: 
\begin{eqnarray}
    \tilde{\pi}_m(w) & = &\mathbb{E}\left[{(r-w)\min\left\{D, \tilde{x}(w) \right\} } \right]=     \mathbb{E}\left[{(r-w)\min\left\{\bb+\A,\bb+\frac{1}{\lambda} \ln\left(\frac{w-k}{c}\right)\right\}}\right] \nonumber \\
     & = & (r-w)\left(\bb+\frac{1}{\lambda}\left(1-\frac{c}{w-k}\right)\right). \label{pimd}
\end{eqnarray}
Hence, the OEM's problem is: $\max_{w \ge c +k}  \  \tilde{\pi}_m(w)$.
By noting that our assumptions
that $r - k - c > c$ and $\frac{1}{\lambda} \ge \bb$ imply $r- k > (\lambda\bb+1)c$, the first-order condition yields
the optimal wholesale price as follows:

\begin{proposition}\label{prop:ws}
In a decentralized supply chain, the optimal wholesale price $\tilde{w}$ satisfies
\begin{equation*}
    \tilde{w} = k+ \sqrt{\frac{(r-k)c}{\bb\lambda+1}},
\end{equation*}
and the corresponding capacity (selected by the supplier) is:
\begin{equation*}
    \tilde{x} \equiv \tilde{x}(\tilde{w}) =\bb+\frac{1}{\lambda} \ln\left(\sqrt{\frac{r-k}{(\bb\lambda+1)c}}\right).
\end{equation*}
\end{proposition}
\noindent
Because $r - k > (\lambda\bb+1)c$, it is easy to check that $\tilde{w}>c + k$ so that the supplier's participation constraint is satisfied. Also, Proposition 2 has the following implications.
First, direct comparison between $\tilde{x}$ and $x^*$ given in Proposition 1 reveals that $\tilde{x} < x^*$.  Hence, when operating under the wholesale price contract $\tilde{w}$ that maximizes the OEM's profit, the supplier will invest capacity level $\tilde{x}$ that is lower than the first best solution $x^*$ for the centralized case.
Second, it is easy to check that both the OEM's optimal wholesale price $\tilde{w}$ and the supplier's ``additional" capacity $\tilde{x} - \bb$ are increasing in the total margin $r-k$ on the end product and the expected additional demand $\frac{1}{\lambda}$, but decreasing in the base demand $\bb$.  
 
By substituting the OEM's optimal wholesale price $\tilde{w}$ and the supplier's capacity $\tilde{x}$ given in Proposition 2 into (\ref{pisd}) and (\ref{pimd}), we get:
\begin{proposition}\label{prop:ws_p}
When the OEM offers its optimal wholesale price $\tilde{w}$ to the supplier in a decentralized supply chain, the supplier's optimal profit $\tilde{\pi}_s$ and the OEM's optimal profit $\tilde{\pi}_m$ satisfy: 
\begin{align*}
    \tilde{\pi}_s &= \left(\sqrt{\frac{(r-k)c}{\bb\lambda+1}}-c\right)\bb + \frac{\sqrt{\frac{(r-k)c}{\bb\lambda+1}}}{\lambda}-\frac{c}{\lambda}-\frac{c}{\lambda}\ln\left(\sqrt{\frac{r-k}{(\bb\lambda+1)c}} \right) > 0, \mbox{and}\\
    \tilde{\pi}_m &=\left(r-k-\sqrt{\frac{(r-k)c}{\bb\lambda+1}}\right) \left(\bb + \frac{1}{\lambda}\left(1-\sqrt{\frac{(\bb\lambda+1)c}{r-k}}\right)\right) > 0.
\end{align*}
Also, $\tilde{\pi}_s+\tilde{\pi}_m < \Pi^*$, where $\Pi^*$ is the optimal profit of the centralized supply chain as stated in Proposition 1. 
\end{proposition}
\noindent
Propositions 2 and 3 reveal that, when operating a decentralized supply chain by using a wholesale price $\tilde{w}$, the supplier will invest at a lower capacity level $\tilde{x}$ (than the first best solution $x^*$) so that the total profit of the decentralized supply chain is lower than that of the centralized supply chain; i.e., $\tilde{\pi}_s+\tilde{\pi}_m < \Pi^*$.

\begin{figure}
    \centering
    \includegraphics[width=0.5\linewidth]{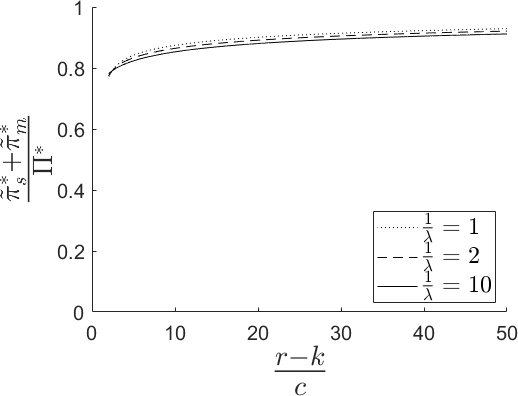}
    \caption{Efficiency of the wholesale price contract $\tilde{w}$.}
    \label{fig:wholesale}
\end{figure}

We now numerically examine the efficiency of the wholesale price contract $\tilde{w}$ (measured in terms of $\frac{\tilde{\pi}_s+\tilde{\pi}_m}{\Pi^*}$) in Figure \ref{fig:wholesale}. Specifically, we set $\bb = 1$, but we vary the ratio between gross margin and the capacity investment cost per unit $\frac{r-k}{c}$ from 2 to 50,\footnote{Because of our assumption $(r-k-c) > c$, $\frac{r-k}{c} > 2$.} and set expected additional demand $\mathbb{E}[A] = \frac{1}{\lambda} = 1$, $2$, $10$.  
\noindent
Observe from Figure \ref{fig:wholesale} that the wholesale price contract $\tilde{w}$ can be rather inefficient when expected additional demand $\frac{1}{\lambda}$ becomes large, which is the case in the high-tech supply chain context we consider.  When $\frac{r-k}{c}$ increases, the inefficiency reduces somewhat, which is mainly due to the fact that the $r-k$ term becomes dominant.

In summary, we find that the wholesale price contract that coordinates the supply chain by setting $ w = r$ is not practical, and the wholesale price contract that optimizes the OEM's profit as stated in Proposition 2 is deemed inefficient (Figure \ref{fig:wholesale}). These observations motivate us to consider ways to ``refine'' the wholesale price contract in order to improve its efficiency in the next section.

\subsection{Augmented wholesale price contracts with lump-sum contingent penalty}\label{sec:penaltycontract}
Recognizing that the traditional wholesale price contract cannot coordinate the supply chain (unless we set $w = r$) and it is inefficient when we set $w = \tilde{w}$ to maximize the OEM's profit, we now consider an ``augmented'' wholesale price contract that can be described as follows: in addition to the wholesale price $w$, the OEM imposes a ``lump-sum contingent'' penalty $\rho$ that the supplier is liable to pay to the OEM when its capacity $x$ is insufficient to fulfill the realized demand $D$; i.e., when $D > x$. Once the wholesale price and penalty are known, the supplier determines its optimal capacity.
We shall show that this augmented wholesale price contract is optimal: it can coordinate the supply chain and it enables the OEM to attain the first best profit $\Pi^*$ as in the centralized system.

To analyze the augmented wholesale price contract $(w, \rho)$ via backward induction, let us first determine the supplier's expected profit. In preparation, let us define an indicator function $1_{\{D > x\}}$ that equals $1$ if $D > x$ and equals $0$, otherwise.  Because $D = \bb + A$, the contingent lump-sum penalty can be expressed as $\rho \cdot 1_{\{\bb+\A > x\}}$. By incorporating this penalty into the supplier's profit given in (\ref{eq:dc_ps}), the supplier's profit for any given augmented wholesale price contract $(w, \rho)$ and any capacity $x$ is $\hat{\pi}_s(x,w,\rho)$, where:
\begin{equation}\label{eq:pis_xwr}
    \hat{\pi}_s(x,w,\rho)= -cx+ \expect{(w-k) \min(\bb+\A, x) - \rho 1_{\{\bb+\A > x\}}}.
\end{equation}
By differentiating $\hat{\pi}_s(x,w,\rho)$ in \eqref{eq:pis_xwr} with respect to $x$ and by considering the first order condition, the optimal supplier's capacity $\hat{x}(w,\rho)$ for any given augmented wholesale price contract  $(w, \rho)$ is:
\begin{equation} \label{xhat}
    \hat{x}(w,\rho) = \bb + \frac{1}{\lambda} \ln \left(\frac{w-k+\rho\lambda}{c}\right).
\end{equation}
Anticipating the supplier's capacity $\hat{x}(w,\rho)$, the OEM's profit is $\hat{\pi}_m(w,\rho)$, where:
\begin{equation}
    \hat{\pi}_m(w,\rho)= \expect{(r-w)\min\{\bb+\A,\hat{x}(w,\rho)\}+ \rho 1_{\{\bb+\A \ge \hat{x}(w,\rho)\}}}
\end{equation}

By comparing the supplier's capacity $\hat{x}(w,\rho)$ against the first best solution $x^*$ given in Proposition 1, we can identify the following conditions for the augmented wholesale price contract $(w, \rho)$ to coordinate the decentralized supply chain so that $\hat{x}(w,\rho) = x^*$:
\begin{proposition}\label{prop:penaltycontract}
Any augmented wholesale price contract $(w, \rho)$ that satisfies $w+\rho \lambda = r$ can enable the OEM to coordinate the decentralized supply chain so that $\hat{x}(w,\rho) = x^*$. However, among all coordinated augmented contracts,
it is optimal for the OEM to set the wholesale price $\hat{w} =k+c+\frac{c}{\bb\lambda+1}\ln\left(\frac{r-k}{c}\right)$ and the contingent lump-sum penalty $\hat{\rho} = \frac{1}{\lambda}\left(r-k-c-\frac{c}{\bb\lambda+1}\ln\left(\frac{r-k}{c}\right) \right) \ge 0$ so that the OEM can extract the entire surplus from the supplier; i.e., $\hat{\pi}_s = 0$, $\hat{\pi}_m=\Pi^*$, and $\hat{\pi}_s + \hat{\pi}_m=\Pi^*$.  
\end{proposition}
\noindent
Proposition \ref{prop:penaltycontract} reveals that there are infinitely many augmented wholesale price contracts  that can coordinate a decentralized supply chain so that the supplier will invest its capacity $\hat{x}(w, \rho) = x^*$.  Also, there exists a coordinating contract $(\hat{w},\hat{\rho})$ that can enable the OEM to extract the entire surplus from the supplier so that her profit is equal to the profit of the entire centralized controlled supply chain. By achieving the first best solution and the highest possible profit, we can conclude that the coordinating contract $(\hat{w},\hat{\rho})$ is the optimal contract for the decentralized supply chain with 100\% contract efficiency (i.e., $\frac{\hat{\pi}_s + \hat{\pi}_m}{\Pi^*} = 1$).

Proposition \ref{prop:penaltycontract} specifies the optimal augmented wholesale price contract for the OEM; however, the implementation of such optimal contract would depend on the underlying business environment.  For instance, the lump-sum penalty $\hat{\rho} = \frac{1}{\lambda}\left(r-k-c-\frac{c}{\bb\lambda+1}\ln\left(\frac{r-k}{c}\right) \right)$ can be too high for a financially constrained supplier to pay so that such an optimal augmented wholesale price contract is not enforceable. This situation can occur in the high-tech sector when the OEM's retail price $r$ is much larger than $c+k$ as illustrated in the following numerical example.

\noindent 
{\bf{A numerical example.}}
Consider the case when $r = 10^7$, $c=10^5$, $k=0$, the base demand $\bb = 50$ and the expected additional demand $\mathbb{E}[A] = \frac{1}{\lambda} = 100$ so that our assumptions $r-k-c > c$ and  $\frac{1}{\lambda} \ge \bb$ hold.  By substituting these parameter values into  Proposition~\ref{prop:penaltycontract}, we get:  $\hat{x} = 510$, $\hat{w} \approx 0.4\cdot 10^6$, and $\hat{\rho} = 959 \cdot 10^6$. 
Now consider a specific realization of the additional demand where $\A = 460$ so that the realized demand  $D = \bb + \A = 510 = \hat{x}$. This demand realization represents the best-case scenario for the supplier, under which he can use its entire capacity to fulfill demand without subjecting to any penalty. The profit for the supplier under this specific demand realization is equal to $156.7 \cdot 10^6$. However, considering all possible demand realizations, because $Prob\{\A > 460\} \approx 0.01$, there is a $1\%$ chance that $D = \bb + A > 510 = \hat{x}$.
When this happens, the supplier is subject to a penalty $\hat{\rho} = 959 \cdot 10^6$ that is six times his best-case profit of $156.7 \cdot 10^6$ and the supplier is unlikely to be able to pay. Therefore, this example illustrates that, when the retail price $r$ is much higher than $k + c$, the optimal augmented wholesale price contract may not be enforceable even though it is optimal for the OEM.

The above numerical example reveals that, when $(r - k -c)$ is very high, the optimal augmented wholesale price contract may not be enforceable because a ``direct'' lump-sum penalty is too high for a financially constrained supplier to pay. 
One could worry that the risk imposed by the supplier is too high in general, since the supplier is held responsible for not meeting uncertain market demands.
However, we need to put this in perspective by considering the bigger picture. Suppose that the supplier delivers a range of relatively inexpensive components to the OEM. For each component the supplier and OEM enter an appropriate contract with contingent penalty. Then obviously, for each single component there may be stockouts that the supplier does not have full control over, but still has to pay for. These stockouts are caused by outside influences beyond the control of the supplier, such as exogenous stochastic demand. However, the supplier still controls how large the risk is. Since the risk of stockouts is low for every component, it is likely that the supplier's actual penalty payments are low compared to his earnings. So overall this is fair.
Another option would be switching to a per-unit penalty. In the following section we will investigate whether this would make the augmented wholesale price contract more applicable in high-tech supply chains.

\subsection{Augmented wholesale price contracts with contingent unit penalty}\label{sec:penaltycontract_unit}
When we consider a per-unit shortfall penalty instead of a lump-sum shortfall penalty, the expected penalty no longer equals  $\rho\mathbb{E}[1_{\{\bb+\A > x\}}]$, but is equal to $\rho_1 \mathbb{E}\left[(b+A-x)^+\right]$ with $\rho_1$ the ``per-unit shortfall penalty''. The supplier's profit function given by \eqref{eq:pis_xwr} is thus adjusted to
\begin{equation}\label{eq: supplier profit unit penalty}
    \hat{\pi}_s(x,w,\rho_1)= -cx+ \expect{(w-k) \min(\bb+\A, x) - \rho_1(b+A-x)^+}
\end{equation}
By considering the first order condition, the optimal supplier's capacity $\hat{x}(w,\rho_1)$ for any given augmented wholesale price contract  $(w, \rho_1)$ is:
\begin{equation} \label{xhat unit}
    \hat{x}(w,\rho_1) = \bb + \frac{1}{\lambda} \ln \left(\frac{w-k+\rho_1\lambda}{c}\right).
\end{equation}
In anticipation of the supplier's capacity $\hat{x}(w,\rho_1)$, the OEM's profit is given by:
\begin{equation}
    \hat{\pi}_m(w,\rho_1)= \expect{(r-w)\min\{\bb+\A,\hat{x}(w,\rho_1)\}+ \rho_1(b+A-\hat{x}(w,\rho_1))^+}
\end{equation}

\begin{proposition}\label{prop: unit penalty}
The augmented wholesale price contract $(w,\rho_1)$ with $\hat{w}=k+c+\frac{c}{\bb\lambda+1}\ln\left(\frac{r-k}{c}\right)$ and $\hat{\rho}_1=r-k-c-\frac{c}{\bb\lambda+1}\ln\left(\frac{r-k}{c}\right)$ coordinates the supply chain while allowing the OEM to extract the entire surplus.
\end{proposition}

Proposition \ref{prop: unit penalty} shows that besides the coordinating augmented contract with lump-sum shortfall penalty $\rho$ presented in Proposition \ref{prop:penaltycontract}, there also exists a coordinating augmented wholesale contract with a per-unit shortfall penalty $\rho_1$ that allows the OEM to capture the entire surplus. Hence, we can conclude that $(\hat{w},\hat{\rho}_1)$ is also an optimal contract for the decentralized supply chain with 100\% contract efficiency. However, returning to the numerical example, in the following we will show that this contract is equally unlikely to be enforceable in a high-tech setting as the lump-sum penalty contract.

\noindent 
{\bf{A numerical example.}} 
Let us revisit the example discussed in \S\ref{sec:penaltycontract}, with $r = 10^7$, $c=10^5$, $k=0$, base demand $\bb = 50$ and expected additional demand $\mathbb{E}[A] = \frac{1}{\lambda} = 100$. By substituting these parameter values into the expressions obtained in Proposition \ref{prop: unit penalty}, we again get $\hat{x}=510$ and $\hat{w}\approx 0.4\cdot 10^6$. We now also have a per-unit penalty $\hat{\rho}_1=9.59\cdot 10^6$, meaning that for every component that the supplier is unable to supply, he incurs a penalty that is nearly 24 times as high as the wholesale price he receives for every unit supplied.  

This numerical example illustrates that when $(r - k -c)$ is very high, the optimal augmented wholesale price contract may not be enforceable, whether a lump-sum penalty is used or a per-unit penalty. Upon discussing with the OEM, we discovered a different form of ``indirect''  penalty that is enforceable by the OEM when it sources components of multiple product generations over time.  We explore this indirect penalty next.

\section{Renewable wholesale price contracts for multiple generations }\label{sec: multi generation}
We now consider the case when the OEM develops multiple generations of a high-tech product.
As explained in \S1, because different generations are based on different technologies, the supplier's capacity is generation-specific:  the extended use of 
the capacity designated for one generation to the next generation is not possible.  
Therefore, to produce components for a new product generation, the supplier needs to invest in new capacity.  Because the capacity is generation-specific, the OEM has the option to work with different suppliers for different generations if the incumbent supplier's performance is unsatisfactory.\footnote{In the high-tech industry, it is common practice for the OEM to work with the incumbent supplier to ensure continuity and smooth transition between product generations unless the supplier's performance is unsatisfactory.}

We learnt from an OEM in the Netherlands that, even though there is an implicit understanding that the OEM would normally renew its contract with the incumbent supplier for the next generation, there is no explicit commitment for contract renewals and there are no explicit conditions for contract non-renewals.  This revelation motivates us to examine a class of wholesale price contracts with ``contingent renewals'': the OEM will renew the contract with the incumbent supplier for the next generation only if the supplier can fulfill the OEM's demand for the current generation.  By specifying the condition for renewal/non-renewal explicitly, the OEM can use contract non-renewal as an ``indirect'' penalty that the OEM can enforce (as opposed to the lump-sum penalty that may not be enforceable). 

In this section, we shall extend our single generation model presented in \S3  to the multi-generation case by incorporating the issue of contingent contract renewals as described above.  Our intent is to examine  the coordinating capability and the efficiency of the contingent wholesale price contract.  
To obtain tractable results, we shall assume that  the demand for each generation $D_t=\bb_t+\A_t$ with $\bb_t=\bb$ for every product generation $t$ and $\A_t$ are i.i.d. exponential random variables with mean $\expect{\A_t} = \frac{1}{\lambda}$. Also, the cost structure is generation-independent so that $r_t=r$, $c_t=c$ and $k_t=k$ for every product generation $t$.\footnote{For the case when the demand and the cost structure are generation-dependent, the analysis is intractable and we shall relegate to future research.}

\subsection{Centralized supply chain with renewals}\label{sec:centralizedmulti}
We first establish a benchmark by considering a centrally controlled supply chain for producing multiple generations with renewals.  Because the demand and the cost structure are generation-independent, the optimal capacity decision is also generation-independent.  Hence, for any capacity $x$ invested for any generation $t$, the profit obtained by the centralized supply chain for this generation is equal to  $\Pi(x)$ as stated in \eqref{eq:sp_sc_p}.  Then, as we consider multi-generation by using a discount factor $0<\delta <1$, the net present value (NPV) of the total supply chain profit with renewals over all generations is:
\begin{equation}
    \Pi^\delta(x) =\sum_{t=0}^\infty \delta^t \Pi(x) = \frac{1}{1-\delta} \Pi(x) = \frac{1}{1-\delta} \cdot \left( -cx + \expect{(r-k)\min\{\bb+\A,x\}} \right).
\end{equation}
By considering the first order condition and by using the fact that $\A\sim Exp(\lambda)$ along with our assumptions
that $r - k - c > c$ and $\frac{1}{\lambda} \ge \bb$, we can apply Proposition \ref{prop:sp_centralized} to show that the optimal capacity for any generation is $x^* = \bb+ \frac{1}{\lambda}\ln\left(\frac{r-k}{c}\right)$. Hence, due to the generation-independent cost structure and identically distributed demand per generation, the same optimization problem is faced for every generation. Therefore, the optimal capacity decision is the same as in the single-generation centralized supply chain. Similarly, the optimal discounted total profit  for the centralized supply chain (associated with the discount factor $\delta$) is denoted by $\Pi^{\delta}$, where 
$\Pi^{\delta} =\frac{1}{1-\delta}\left((r-k-c)\left(\bb+\frac{1}{\lambda}\right) - \frac{c}{\lambda}\ln\left(\frac{r-k}{c}\right)\right) \ge 0$.

\subsection{Wholesale price contracts with exogenous renewal probability}\label{sec:multiExogenous}

We now consider a decentralized supply chain in which an OEM establishes a supply contract with a focal supplier.  
Similar to the decentralized case for a single generation, the OEM determines the wholesale price, after which the supplier decides on capacity. In addition, we have the possibility for the OEM to renew the contract with the incumbent supplier. The sequence of decision making and actions can thus be summarized as:
\begin{enumerate}
    \item OEM determines wholesale price
    \item Supplier invests in production capacity
    \item Demand is realized and supplier produces required components up to the maximum capacity
    \item OEM decides whether or not to renew contract with supplier
\end{enumerate}

In \S4.3, we shall analyze the contingent wholesale price contract under which the renewal for the next generation depends on whether the supplier's capacity can satisfy the OEM's demand for the current generation.  To explicate our analysis, let us first consider a base case in which contract renewal is based on an ``exogenously'' given real-valued probability $R$ that is independent of the supplier's capacity decision $x$. Once a contract is not renewed, we assume that the supplier will never be allowed to work with the OEM in the future. This assumption implies that the number of generations $Y$ that the incumbent supplier can work with the OEM will follow a geometric distribution so that $Y \sim Geom(1-R)$ and $Prob\{Y = t\} = R^{t-1}(1-R)$ for $ t = 1, 2, \dots$. Also, observe that the supplier's capacity decision $x$ is generation-independent because the distribution of demand $D$ and all cost parameters are generation-independent. Consequently, the expected profit for the supplier in each generation $t$ is $\tilde{\pi}_s(x ,w)$ as stated in \eqref{eq:dc_ps}.   By combining these observations, the NPV of the supplier's expected profit over $Y$ product generations can be expressed as:
\begin{align}\label{eq:mp_s_p0}
    \pi_s^\delta(x)=\expect{\sum_{t=1}^{Y}\delta^{t-1} \tilde{\pi}_s(x,w)}=  \sum_{t=1}^{\infty} P(Y \ge t) \delta^{t-1}\tilde{\pi}_s(x,w)   =   \sum_{t=1}^{\infty} (R\delta)^{t-1} \tilde{\pi}_s(x,w) = \frac{\tilde{\pi}_s(x,w)}{1-\delta R}.
\end{align}

\noindent
It follows from \eqref{eq:mp_s_p0} that the term $\frac{1}{1-\delta R}$ is independent of $x$, so we can conclude that, for any given $R$, the NPV of the supplier's profit is maximized when he maximizes his single-period profits $\tilde{\pi}_s(x,w)$. Combining this observation with our analysis presented in \S3.2, we can conclude that, when the contract renewal probability $R$ is exogenously given, it is optimal for the supplier to set its capacity according to $\tilde{x}(w) = \bb+\frac{1}{\lambda}\ln(\frac{w-k}{c})$ as stated in \eqref{eq:sp_dc_x}, where $\tilde{x}(w)$ is independent of $R$. Since the supplier's capacity investment does not affect the renewal decision, there is no incentive for the supplier to invest in additional capacity and the optimal capacity remains the same as in the single-generation case.
Also, we can use the same approach as presented in \S3.1 (Corollary 1) to show that a wholesale price contract with exogenously given renewal probability $R$ can coordinate the supply chain (i.e., $\tilde{x}(w) = x^*$) only when we set $w = r$, which the OEM will not oblige.  

Instead of coordinating the supply chain, the OEM can seek to maximize her own profit.  Because the demand $D$ and all cost parameters are generation-independent, we can use the same approach to show that the optimal wholesale price $\tilde{w}$ is given in Proposition \ref{prop:ws}. Thus the corresponding contract is inefficient (cf. Proposition \ref{prop:ws_p}). 

In cognizant of the shortcomings of the wholesale price contracts with exogenous renewal probability $R$, we now examine the coordinating capability and the efficiency of the wholesale price contracts with renewal probability $R(x)$ that is ``endogenously'' determined by the supplier.  

\subsection{Wholesale price contracts with endogenous renewal probability}\label{sec:performconting}

We now extend our analysis presented in the previous section to the case when the contract will be renewed for the next generation only if the supplier's capacity $x$ can meet with the OEM's demand $D$ for the current generation. This way, the renewal probability $R(x)$ is now ``endogenously'' dependent on the capacity $x$ to be selected by the supplier.  Based on our assumptions that $D = \bb +A$ and $A \sim Exp(\lambda)$ are generation-independent, it is easy to check that, for any given supplier capacity $x$, the renewal probability $R(x) = Prob\{ D \le x \} = Prob\{ A \le  (x - \bb)\} = 1-e^{-\lambda(x-\bb)}$ for $x\geq \bb$, and $R(x)=0$; otherwise.  Hence, for any given supplier capacity $x$, we can use the same approach as presented in \S4.2 to show that the number of generations $Y(x)$ that the incumbent supplier can work with the OEM will follow a geometric distribution so that $Y \sim Geom(1-R(x))$, where the renewal probability $R(x)$ is defined above.  As before, because the supplier's capacity decision $x$ is generation-independent (because the demand $D$ and all cost parameters are generation-independent), the expected profit for the supplier in each generation $t$ is $\tilde{\pi}_s(x ,w)$ as stated in \eqref{eq:dc_ps}. These observations imply that the NPV of the supplier's expected profit over $Y(x)$ product generations can be expressed as: 
\begin{align}\label{eq:mp_s_p}
    \tilde{\pi}_s^\delta(x)=\expect{\sum_{t=1}^{Y(x)}\delta^{t-1} \tilde{\pi}_s(x,w)}=\sum_{t=1}^{\infty}P( Y(x)\ge t)\delta^{t-1} \tilde{\pi}_s(x,w)=\sum_{t=1}^{\infty}(R(x)\delta)^{t-1} \tilde{\pi}_s(x,w)=\frac{\tilde{\pi}_s(x,w)}{1-\delta R(x)},
\end{align}
where the renewal probability $R(x) = 1-e^{-\lambda(x-\bb)}$ for $x\geq \bb$, and $R(x)=0$; otherwise. 

By considering the first-order condition and the assumption that $w \ge c+k$, we can determine the supplier's optimal capacity $\tilde{x}(w, \delta)$ as follows: 

\begin{proposition}\label{prop:sup_opt}
For any given wholesale price contract $w$ (that has $w \ge c + k$) with endogenous renewal probability $R(x)$, 
the supplier's optimal capacity satisfies: 
\begin{equation}\label{eq:sup_opt}
\tilde{x}(w,\delta)=\bb+\frac{1}{\lambda} \ln\left(\frac{\delta}{1-\delta}W\left( \frac{1-\delta}{\delta} e^{ \frac{w-k}{\delta c}-1 + \frac{\bb\lambda(w-k-c)}{c}} \right)\right)
\end{equation} 
where $W(\cdot)$ is known as the Lambert $W$ function.\footnote{The Lambert $W(y)$ function is the inverse function of $f(y)=ye^y$.} 
\end{proposition}
\noindent
Analogous to the supplier's optimal capacity $\tilde{x}(w)$ given in \eqref{eq:sp_dc_x} for the single-generation case, the supplier's optimal capacity  $\tilde{x}(w,\delta)$ is equal to the base demand $\bb$ plus some additional capacity to cover the uncertain additional demand $A$. Unlike $\tilde{x}(w)$, observe from \eqref{eq:sup_opt} that the additional capacity $\frac{1}{\lambda} \ln\left(\frac{\delta}{1-\delta}W\left( \frac{1-\delta}{\delta} e^{ \frac{w-k}{\delta c}-1 + \frac{\bb\lambda(w-k-c)}{c}} \right)\right)$
created by  $\tilde{x}(w,\delta)$ under the ``contingent'' wholesale price contract depends on the base demand $\bb$.  The reason lies in the fact that, when $\bb$ increases, the supplier values contract renewals more because it can obtain a higher profit through the base demand. Consequently, as the base demand becomes bigger, the supplier has stronger incentive to invest in more capacity to increase its renewal probability. Besides the impact of the base demand $\bb$, it is easy to show that as the discount factor $\delta$ increases, the supplier will value contract renewals more; hence, the supplier will increase its capacity as $\delta$ increases so that $\tilde{x}(w,\delta)$ is increasing in $\delta \in (0, 1)$. Thus high-tech manufacturers can benefit from seeking out suppliers that have a focus on the long-term profit. 

\begin{figure}
    \centering
    \includegraphics[width=0.5\textwidth]{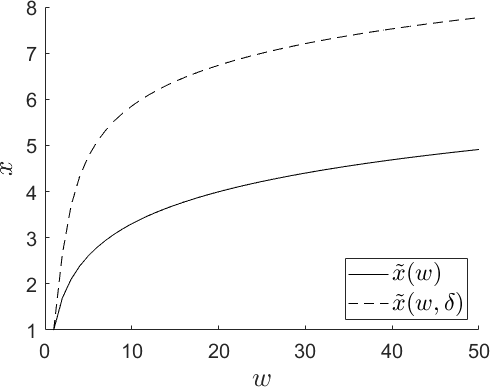}
    \caption{Supplier's optimal capacity for exogenous and endogenous renewal probabilities for different wholesale prices, with $b=1$, $\lambda=1$, $\delta=0.9$, $c=1$ and $k=0$.}
    \label{fig: capacity comparison}
\end{figure}

Next, through direct comparison between $\tilde{x}(w)$ given in \eqref{eq:sp_dc_x} for the single-generation case (which corresponds to the multi-generation case with exogenous renewal probability $R$ and yet $\tilde{x}(w)$ is independent of $R$ as explained in \S4.2) and $\tilde{x}(w,\delta)$ given in \eqref{eq:sup_opt} for the multi-generation case with endogenous renewal probability $R(\tilde{x}(w,\delta))$, we get: 
\begin{proposition}\label{prop:sl}
For any given wholesale price contract $w$ (that has $ w \ge c + k$), the supplier's optimal capacity  $\tilde{x}(w,\delta)$ corresponding to endogenous renewal probability $R(\tilde{x}(w,\delta))$ is higher than its optimal capacity $\tilde{x}(w)$ corresponding to any exogenous renewal probability $R$ (i.e., $\tilde{x}(w,\delta)>\tilde{x}(w)$). 
\end{proposition}
\noindent
While the supplier's capacity (i.e., both $\tilde{x}(w,\delta)$ and $\tilde{x}(w)$) is increasing in the wholesale price $w$, the above proposition suggests that the stipulated condition for contract renewal provides an incentive for the supplier to invest more capacity.
This is also illustrated in Figure \ref{fig: capacity comparison}, where we show the supplier's optimal capacity both in case of an endogenous renewal probability ($\tilde{x}(w,\delta)$) and in case of an exogenous renewal probability ($\tilde{x}(w)$) for different values of wholesale price $w$, with $b=1$, $\lambda=1$, $\delta=0.9$, $c=1$ and $k=0$. We observe that the optimal capacity under the contract with endogenous renewal probability is indeed considerably higher than in case of an exogenous renewal probability, for the same wholesale price.

\subsubsection{Supply chain coordination} \label{sec: endogenous renewal - coordinating wholesale}
Recall from \S4.2 that the wholesale price contract with exogenous renewal probability $R$ can coordinate the supply chain (i.e., $\tilde{x}(w) = x^*$) only when the OEM sets $w = r$, which the OEM will not oblige because of zero profit.  We now examine whether the use of endogenous renewal probability $R(x)$ would enable the OEM to coordinate the supply chain.  By considering the supplier's optimal capacity $\tilde{x}(w,\delta)$ given in \eqref{eq:sup_opt} and $x^*$ given in Proposition 1, we get: 
\begin{proposition}\label{prop:mp_dc_c_w}
Suppose the OEM sets its wholesale price for each generation at $w^\delta$, where:
\begin{equation*}
    w^\delta = k+\frac{\delta c \left(1+\bb\lambda+\ln\left(\frac{r-k}{c}\right)\right)+(1-\delta)(r-k)}{1+\delta \bb\lambda} < r.
\end{equation*}
Then corresponding contingent wholesale price contract with endogenous renewal probability $ R(\tilde{x}(w^\delta,\delta))$ can coordinate the supply chain so that $\tilde{x}(w^\delta,\delta)=\bb+\frac{1}{\lambda}\ln\left(\frac{r-k}{c}\right)=x^*$. 
\end{proposition}
\noindent
The above proposition reveals that, by imposing an explicit condition for contract renewal with the incumbent supplier, the OEM can leverage the indirect penalty associated the non-renewal to entice the supplier to select its capacity according to the first best solution $x^*$ by offering the wholesale price $w^\delta$. Since $w^{\delta}<r$, supply chain coordination is achieved at a strictly lower wholesale price compared to the wholesale price of $r$ that is required for coordination when renewal is exogenous (cf \S\ref{sec:multiExogenous}). 

\begin{figure}
    \centering
    \includegraphics[width=0.5\textwidth]{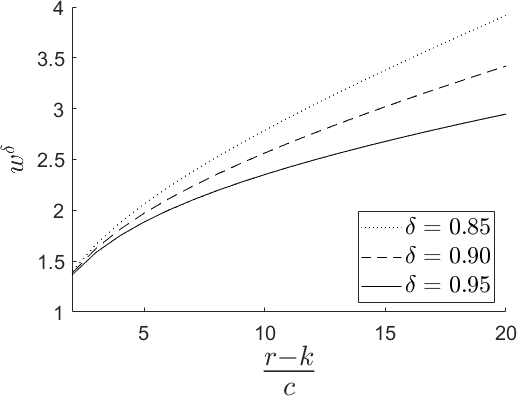}
    \caption{Coordinating wholesale price as a function of $\frac{r-k}{c}$ for different values of $\delta$}
    \label{fig: coordinating w}
\end{figure}

In Figure \ref{fig: coordinating w} we analyze numerically how the coordinating wholesale price changes with the margin on the end-product for different discount factors, where $c=1$, $k=0$, $b=1$ and $\lambda=1$. We observe that the coordinating wholesale price that the OEM pays to the supplier increases with $\frac{r-k}{c}$. When the margin on the end product is larger, the supplier will require a higher wholesale price to build the first-best capacity $x^*$. However, when the valuation of future profits by the supplier (represented by $\delta$) is larger, the increase in the coordinating wholesale price as $\frac{r-k}{c}$ increases becomes smaller as the supplier is more inclined to invest in future profits.

\subsubsection{Supplier Surplus Extraction} \label{sec: endogenous renewal - surplus extraction}
\begin{figure}
    \centering
    \includegraphics[width=\linewidth]{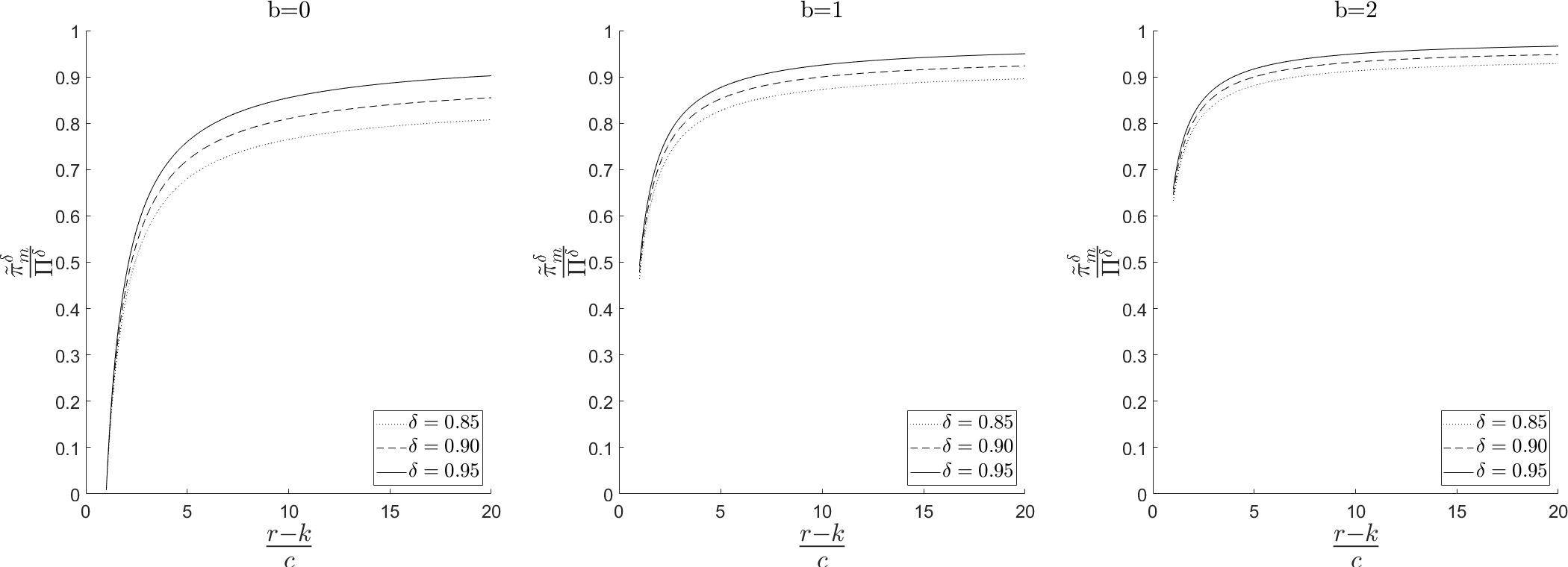}
    \caption{Fraction of the NPV captured by the OEM ($\tilde{\pi}^\delta_m/\Pi^\delta$) under the coordinating contingent wholesale price contract $(w^\delta)$, for different values of the base demand $\bb$.}
    \label{fig:distr_profits}
\end{figure}

Proposition \ref{prop:mp_dc_c_w} shows that the contingent wholesale price contract can coordinate the supply chain by offering a wholesale price $w^\delta < r$.  
We now examine whether such coordinating contract can enable the OEM to extract the entire surplus from the supplier so that the corresponding contract is optimal.  In preparation, let us first compute the NPV of the OEM's profit over all generations.  For any given supplier capacity $x$ and wholesale price $w$, the OEM will earn a profit $\tilde{\pi}_{m}(x,w) = \mathbb{E}[{(r-w)\min\{\bb + A, x \} }]$ for each generation regardless of which supplier its works with.  In other words, even though the OEM may work with different suppliers upon contract non-renewals, the OEM's profit for each generation stays the same and it is independent of the contract renewal probability $R(x)$.  Because of the assumptions that the demand $D = \bb + A$ and all cost parameters are generation-independent, the NPV of the OEM's profit over all generations can be expressed as:  
\begin{equation}\label{eq:pi_m}
    \tilde{\pi}_{m}^\delta(x,w)= \sum_{t=0}^\infty \delta^t \tilde{\pi}_{m}(x,w) =\frac{\tilde{\pi}_{m}(x,w)}{1-\delta}
\end{equation}
Now, the OEM offers a wholesale price $w^\delta (< r)$ as stated in Proposition \ref{prop:mp_dc_c_w} (along with the contingent renewal condition) to entice the supplier to set $\tilde{x}(w^\delta,\delta) =x^* = \bb+\frac{1}{\lambda}\ln\left(\frac{r-k}{c}\right)$.  Hence, we can use the fact that $A \sim Exp(\lambda)$ to determine the NPV of the OEM's profit as $\tilde{\pi}_{m}^{\delta}
=  \frac{1}{1-\delta}\left(r-w^\delta\right)\left(\bb+\frac{1}{\lambda}(1-\frac{c}{r-k})\right) > 0$. By comparing $\tilde{\pi}_{m}^{\delta}$ against the optimal NPV of the centralized supply chain $\Pi^\delta$ as defined in \S\ref{sec:centralizedmulti}, we get:
\begin{proposition}\label{prop:NPVcapture}
Under the coordinating contract ($w^\delta$) with endogenous renewal probability, the OEM cannot extract the entire surplus from the supplier: the fraction of the NPV of the total supply chain profit captured by the OEM  $\frac{\tilde{\pi}_{m}^{\delta}}{\Pi^\delta} < 1$. 
\end{proposition}
\noindent
By considering different values of $\bb$ and $\delta$, for $\lambda=1$, Figure \ref{fig:distr_profits} depicts the fraction of the NPV of the total supply chain profit captured by the OEM (given by $\frac{\tilde{\pi}_{m}^{\delta}}{\Pi^\delta}$) under the coordinating contract $(w^\delta)$ with endogenous renewal probability as a function of the margin on the end product. Figure \ref{fig:distr_profits} verifies that the fraction is strictly below 1.  Also, from these figures, we notice that the OEM can capture a larger proportion of the NPV of the total supply chain profit when $\delta$ is large. This is due to the fact that a supplier that has a high valuation of future profits requires less incentive to invest sufficient capacity. Similarly, the supplier is more willing to invest when the base demand, resulting in a certain profit, is higher, which results in a higher fraction of the profit for the OEM. Furthermore, the fraction of the NPV captured by the OEM is higher when his added value and thus $\frac{r-k}{c}$ is large.

Next, to gain analytic insights into whether the coordinating contract with endogenous renewal probability is suitable in the high-tech setting, we investigate the fraction of NPV captured by the OEM for the following special case of our general model that is very relevant in the high-tech setting, as explained in \S\ref{sec:introduction}. \\

\noindent
{\bf{Special Case 1:  When the ratio $\frac{r-k}{c}$ is very large (or the margin $r-k-c$ is very large).}}  The following proposition characterizes the fraction of the optimal NPV captured by the OEM ($\tilde{\pi}^\delta_m/\Pi^\delta$) as defined in Proposition \ref{prop:NPVcapture} when $\frac{r-k}{c}\rightarrow\infty$. This limiting case is of interest because it corresponds to the situation when the penalty under the augmented wholesale price contract is exorbitant to be enforceable as explained in \S\ref{sec:penaltycontract}.\footnote{The extreme case reflects situations where the supplier delivers a crucial part of the system developed by the OEM, with a value that is much lower than the selling price of the product, as is typical in high-tech manufacturing: $r-k$ corresponds to the gross margin when selling the product (e.g. an aircraft, a wafer-stepper), while $c$ corresponds to the costs of capacity for producing a component of that product (e.g. a wing section, a wafer handler).} 

\begin{proposition}\label{prop:distr_profits}
Suppose the OEM offers the coordinating contract ($w^\delta$) with endogenous renewal probability. Then, when $\frac{r-k}{c}\rightarrow\infty$, the fraction of the NPV of the total supply chain profit captured by the OEM 
$\tilde{\pi}^\delta_m/\Pi^\delta \rightarrow  \frac{\delta \bb \lambda+\delta}{\delta \bb\lambda+1}$.
\end{proposition}
\noindent
The above proposition has the following implications. First, when the base demand $\bb = 0$, the limit of the fraction $\tilde{\pi}^\delta_m/\Pi^\delta$ is equal to $\delta$. Hence, when the base demand is low and when the discount factor is high, the OEM can extract the bulk of the surplus from the supplier by adopting the coordinating contingent wholesale price contract with endogenous renewal probability. Second, when the base demand $\bb$ is substantial in comparison to $\mathbb{E}[A]=1/\lambda$, the limit of the fraction $\tilde{\pi}^\delta_m/\Pi^\delta$ will approach $1$. Hence, when a large portion of total demand is obtained through advance orders, the OEM can extract almost the entire surplus from the supplier. This means that having a strong market position, demonstrated by many advance orders, also gives the OEM a strong position vis-à-vis her supplier. 
These two observations enable us to characterize the business environment (i.e., when $\frac{r-k}{c}$ is high, base demand $\bb$ is substantial compared to additional demand $\mathbb{E}[\A]$, or the discount factor $\delta$ is high) in which the coordinating wholesale price contract with contingent renewal can enable the OEM to extract the bulk of the surplus from the supplier so that this contract is close to optimal.

\noindent
{\bf{Special Case 2:  When the ratio $\frac{r-k}{c}$ is close to but strictly greater than 1 (or the margin $r-k-c$ is very small).}} This situation occurs when the supplier's capacity cost $k$ and unit cost $c$ are high, such that $\frac{r-k}{c}\rightarrow 1^+$.\footnote{$\frac{r-k}{c}\rightarrow 1^+$ denotes  $\frac{r-k}{c}$ approaches $1$ from the right so that $\frac{r-k}{c} > 1$.} By considering the case when $\delta \rightarrow 1$ we can compare the coordinating wholesale price $w^\delta$ given in Proposition 7 for the multi-generation case associated with the contingent renewal contract and the coordinating contingent penalty $\hat{\rho}$ given in Proposition 4 for the single generation case as examined in \S3.3.  

\begin{proposition}\label{prop:distr_profits_small_y}  When $\frac{r-k}{c}\rightarrow 1^+$ and $\delta \rightarrow 1$, the coordinating wholesale price $w^\delta$ given in Proposition 7 for the multi-generation case satisfies: $w^\delta \approx r$.  However, the coordinating contingent penalty $\hat{\rho}$ given in Proposition 4 for the single generation case satisfies: $\hat{\rho} \approx 0$.
\end{proposition}

The above proposition has the following implications: when the margin $r-k-c$ is very small and $\delta$ is close to 1, the contingent penalty $\hat{\rho}$ is very small so that the augmented wholesale price contract is easily enforceable.  Also, as revealed in Proposition 4, it allows the OEM to capture the entire supply chain profit. However, the wholesale price $w^\delta$ is close to $r$ so that the contingent renewal contract renders the OEM essentially profitless.  Therefore, when the margin $r-k-c$ is very small and $\delta$ is close to $1$, the OEM is better off to treat each generation separately by adopting the contingent penalty contract instead of the contingent renewal contract for multiple generations.  
  
Based on the analysis of these two special cases, we can make the following conclusions.  First, the long-term supply contract with contingent renewals is effective for the OEM when the margin $r - k -c$ is very large.  This is because, when the margin $r - k -c$ is very large,  the contingent penalty $\hat{\rho}$ is exorbitant so that the augmented contract with contingent penalty is deemed impractical. However, the contingent renewal contract performs well because it enables the OEM to obtain almost the entire supply chain profit.  Second, the short-term contingent penalty contract is more efficient for the OEM when the margin $r - k -c$ is very small and $\delta$ is close to 1. This is because, in this case, the coordinating contingent wholesale price $w^\delta \approx r$, leaving very little profit for the OEM under the contingent renewal contract.  However, the contingent penalty $\hat{\rho}$ is very small so that the augmented contract with contingent penalty can be easily enforced and yet it enables the OEM to obtain the entire supply chain profit.

\subsubsection{Duration of collaboration}\label{sec: duration}
Now that we have established the importance of long-term collaborations in high-tech supply chains, the question can be asked how long these collaborations will last. Since the duration of the collaboration $Y$ is distributed geometrically with parameter $1-R(x)$, where $R(x)$ is the renewal probability, the expected duration of the collaboration equals $\mathbb{E}[Y(x)]=\frac{1}{1-R(x)}=\frac{1}{e^{-\lambda(x-b)}}$. Proposition \ref{prop: duration collaboration coordination} shows that the duraction of the collaboration when the OEM sets coordinating wholesale price $w^\delta$, to induce the supplier to set capacity $x^*=b+\frac{1}{\lambda}\ln\left(\frac{r-k}{c}\right)$, is equal to $\frac{r-k}{c}$. This means that the higher the value of the end-product, the longer the collaboration lasts. Also, due to our assumption that $r-k-c > c$, we get:

\begin{proposition}\label{prop: duration collaboration coordination}
Under the coordinating wholesale price $w^\delta$, the expected duration of the collaboration is $\mathbb{E}[Y(x)]=\frac{1}{1-R(x)}=\frac{r-k}{c} > 2$.
\end{proposition}

\subsubsection{Optimal wholesale price renewal contract}\label{sec: optimal wholesale renewal}
Now that we have shown that a renewal contract with endogenous renewal probability can coordinate the supply chain while yielding a positive profit to both parties, the question remains whether the manufacturer has incentive to set this coordinating wholesale price. Since, optimization of the wholesale price under an endogenous renewal probability is intractable, we answer this question based on numerical experiments. Since the OEM's profit function is concave, the optimal wholesale price can be determined numerically using Golden-section search.

We consider several instances with $c=1$ and take $r-k\in\{10,50,100\}$, $B\in\{0,1\}$, $\lambda\in\{0.1,0.5,1\}$ and $\delta\in\{0.85,0.9,0.95\}$. This gives in total $3^3\cdot 2=54$ instances. For every instance we calculate the profit for the OEM under the contingent renewal contract with endogenous renewal probability for both the coordinating and the OEM's optimal wholesale price and determine the percentage difference. In addition, we determine the expected number of generations that the collaboration will last, both under the optimal and the coordinating wholesale price. 
The summarizing statistics of the full factorial experiment are given in Table \ref{tab: difference profit coordinated}. 

We can observe that in all instances the OEM looses some money by setting the coordinating wholesale price instead of optimizing the wholesale price. When we first consider the effect of $r-k$, we observe that the profit lost by coordinating the supply chain reduces as the value of the end product increases. Furthermore, we observe that $\delta$ has a large effect on the difference in profit. When $\delta$ increases, the difference in profit reduces considerably. The same holds for base demand $b$. The parameter values for which the difference in profit between choosing the optimal and coordinating wholesale price is smallest thus correspond to the case for which the coordinating contract is most suitable, namely a high value end-product and high valuation of future profits.
When we additionally consider the duration of the collaboration, we observe that under both the optimal and coordinating wholesale price the collaboration is expected to span multiple generations. Furthermore, we observe that collaboration lasts considerably longer under the coordinating wholesale price than under the optimal wholesale price, especially for high-valued end-products. 

\begin{table}[]
    \centering
    \caption{Results full factorial experiment on the difference in OEM's profit per period between using the coordinating and optimal wholesale price}
    \label{tab: difference profit coordinated}
    \begin{tabular}{lrrrrrrrr}
    \hline
    & & \multicolumn{2}{c}{ Average profit OEM per period} & \multicolumn{3}{c}{Difference profit (\%)} & \multicolumn{2}{c}{ Average duration relation} \\
    & & Optimal & Coordinated & Average & Max & Min &  Optimal & Coordinated\\
    \hline
$r-k$	&	5	&	9.98	&	9.27	&	6.81	&	11.89	&	1.71  & 3.21   &	5.00	\\
	&	10	&	29.23	&	27.81	&	4.67	&	8.97	&	0.96    &   5.79    &   10.00	\\
	&	20	&	71.91	&	68.55	&	4.46	&	8.78	&	0.90	&   9.86    &   20.00\\
\hline
$\delta$	&	0.85	&	36.36	&	33.46	&	8.15	&	11.89	&	4.07    &   5.55    &   11.67	\\
	&	0.9	&	36.99	&	35.21	&	5.16	&	8.40	&	2.34	&   6.16    &   11.67\\
	&	0.95	&	37.76	&	36.96	&	2.62	&	5.28	&	0.90    &   7.17    &   11.67	\\
\hline
$b$	&	0	&	31.89	&	30.09	&	6.52	&	11.89	&	2.08	&   6.32    &   11.67 \\
	&	1	&	42.19	&	40.33	&	4.11	&	10.37	&	0.90    &   6.27    &   11.67	\\
\hline
$\lambda$	&	0.1	&	78.69	&	74.49	&	6.12	&	11.89	&	1.84	&   6.31    &   11.67\\
	&	0.5	&	19.87	&	19.01	&	5.19	&	11.89	&	1.26	&   6.29    &   11.67\\
	&	1	&	12.56	&	12.12	&	4.63	&	11.89	&	0.90	&   6.28    &   11.67\\
\hline
Overall	&		&	12.56	&	12.12	&	5.31	&	11.89	&	0.90	&   6.29    &   11.67\\
\hline
    \end{tabular}
\end{table}

\section{Extensions}\label{sec: extensions}
This section analyzes the effect of (1) a supplier's reservation profit and (2) more general demand distributions on the effectiveness of the different supply contracts and the corresponding division of profit between the supplier and the OEM.

\subsection{Supplier's reservation profit}\label{sec: extension - reservation profit}
Until now we have assumed that the supplier will engage as long as the expected profit is non-negative. However, it is likely that a supplier will request a positive profit to justify his efforts. Since the contingent renewal contract already guarantees a positive profit for both parties, we investigate in this section how including a reservation profit, denoted by $\reservationprofit$, affects our analysis of the single-generation contracts. 

\subsubsection{Single-generation wholesale price contract}
Since the supplier's reservation profit does not affect the policy parameters, the supplier's profit function as given in \eqref{eq:dc_ps} remains the same. In Lemma \ref{lem:helperlemma}, we analyzed for which values of wholesale price $w$ the supplier's expected profit is non-negative. Analogously, Lemma \ref{lem: reservation profit} gives the minimum wholesale price for which the supplier attains the reservation profit.
\begin{lemma}\label{lem: reservation profit}
The supplier's profit $\tilde{\pi}_s(w)=(w-c-k)\left(b+\frac{1}{\lambda}\right)-\frac{c}{\lambda}\ln\left(\frac{w-k}{c}\right)\geq\reservationprofit$ if and only if $w\geq k-\frac{c}{b\lambda+1}W\left(-(b\lambda+1)e^{-(b\lambda+1+\reservationprofit\frac{\lambda}{c}}\right)$.
\end{lemma}
This means that the supplier will engage in the wholesale price contract proposed in Proposition \ref{prop:ws} provided that the condition given in Lemma \ref{lem: reservation profit} is satisfied. If the OEM's optimal wholesale price does not satisfy this condition, the OEM will need to pay a higher than optimal wholesale price to the supplier, leaving the OEM with lower profits.

\subsubsection{Augmented wholesale price contract with lump-sum penalty}
Under the OEM's optimal augemented wholesale price contract with lump-sum penalty, which was proposed Proposition \ref{prop:penaltycontract} in \S\ref{sec:penaltycontract}, the OEM was able to capture the entire supply chain profit. When the supplier has a reservation profit $\reservationprofit>0$, this is no longer possible. In this case, it will be optimal for the OEM to determine the policy parameters $w$ and $\rho$ for which the supplier builds the first best capacity $x^*$ and the supplier's expected profit is exactly equal to the reservation profit. The details of this optimal policy are given in Proposition \ref{prop:penaltycontract reservation profit}.

\begin{proposition}\label{prop:penaltycontract reservation profit}
Any augmented wholesale price contract $(w, \rho)$ that satisfies $w+\rho \lambda = r$ can enable the OEM to coordinate the decentralized supply chain so that $\hat{x}(w,\rho) = x^*$. However, given the supplier's reservation profit $\reservationprofit$,
it is optimal for the OEM to set the wholesale price $\hat{w} =k+c+\frac{c}{\bb\lambda+1}\ln\left(\frac{r-k}{c}\right)-\reservationprofit\frac{\lambda}{\bb\lambda+1}$ and the contingent lump-sum penalty $\hat{\rho} = \frac{1}{\lambda}\left(r-k-c-\frac{c}{\bb\lambda+1}\ln\left(\frac{r-k}{c}\right) -\reservationprofit\frac{\lambda}{\bb\lambda+1}\right)$ such that $\hat{\pi}_s = \reservationprofit$, $\hat{\pi}_m=\Pi^*-\reservationprofit$, and $\hat{\pi}_s + \hat{\pi}_m=\Pi^*$. Furthermore, $\hat{\rho}\geq 0$ for $Z\leq \Pi^*$.
\end{proposition}

The optimal augmented wholesale price contract with unit penalty can be analyzed accordingly.

\subsection{Erlang-$\erlangshape$ distributed demand}\label{sec: erlang demand}
In the previous sections, we assumed demand to consist of a base demand corresponding to advance orders and an uncertain exponentially distributed part to be able to derive theoretical and numerical results that are relevant in a high-tech setting. We will now investigate numerically whether our findings extend to more general distributions for the uncertain demand part. 
To maintain the characteristics of our problem and the relevance of our results for high-tech industries, we need to focus on demand distributions that only allow for positive demands and have a long tail, since there is a small but positive probability of very large demand. Considering demand distributions with a finite upper bound on the domain leads to a fundamentally different problem, as it can be guaranteed that the supplier has sufficient capacity when the capacity decision equals the maximum demand. Specifically, we will consider $A\sim Erlang(\lambda,\erlangshape)$ such that demand consists of a fixed part and an Erlang distributed part for different values of $\erlangshape$, yielding different coefficients of variation.
Due to the subsequent optimization of capacity and wholesale price with the renewal probability endogenously determined by these values, analytical analysis breaks down and we resort to numerical analysis.

Using the expected sales defined in Lemma \ref{lem: extension - exp sales}, we can express the supplier's and OEM's expected profit functions that we analyze numerically. We use a Golden-section search procedure to find the optimal capacity for a given wholesale price and bisection search to find the coordinating wholesale price.

\begin{lemma}\label{lem: extension - exp sales}
For $D\sim Erlang(\lambda,\erlangshape)$ and capacity $x$, the expected sales are given by:
\begin{equation}
    \mathbb{E}\left[ \min\{D,x\} \right]=b+\mathbb{E}\left[ \min\{A,x-b\} \right]=b+\frac{\gamma(k+1,\lambda (x-b))}{\lambda(\erlangshape-1)!} + (x-b) \mathbb{P}(A>x-b)
\end{equation}
where $\gamma(s,x)=\int_0^x t^{s-1}e^{-t} dt$ is the lower incomplete gamma function.
\end{lemma}

\subsubsection{Penalty contract}
It can be verified numerically that also for Erlang-distributed demand, an augmented wholesale price  contract with lump-sum penalty may not be enforceable in high-tech supply chains. For example, let us return to the numerical example introduced in \S\ref{sec:penaltycontract}, with $r=10^7$, $c=10^5$ and $k=0$ and assume that $b=50$ and $A\sim Erlang\left(\frac{3}{100},3\right)$ such that again $\mathbb{E}[D]=150$. We find that under the OEM's optimal contract $\hat{w}\approx 0.14\cdot10^6$ and $\hat{\rho}\approx25.2\cdot10^6$. The supplier's optimal capacity then equals $\hat{x}\approx191$. When realized demand is exactly equal to $\hat{x}$, the supplier can use its entire capacity and earn $6.7\cdot10^6$, without subjecting to the penalty. However, when realized demand does exceed the available capacity, the supplier is subject to a penalty $\hat{\rho}\approx 25.2\cdot10^6$ that is 3.7 times as high as his best-case profit. Hence, also for $A\sim Erlang\left(\frac{3}{100},3\right)$, the optimal augmented wholesale price contract may not be enforceable when $r$ is much larger than $c+k$.

\subsubsection{Wholesale price contract with endogenous renewal probability}
\begin{figure}
    \centering
    \includegraphics[width=\textwidth]{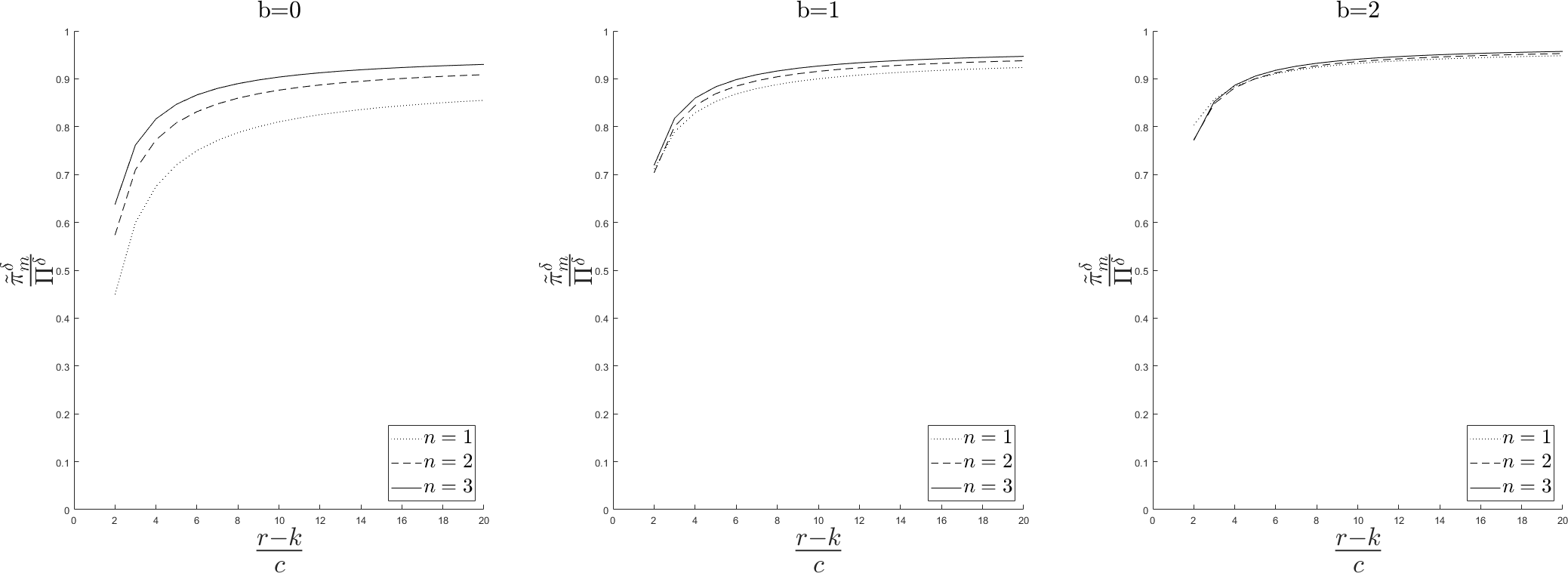}
    \caption{Fraction of the NPV captured by the OEM ($\tilde{\pi}^\delta_m/\Pi^\delta$) under the coordinating contingent wholesale price contract with Erlang-$\erlangshape$ demand.}
    \label{fig: division profits erlang}
\end{figure}

Now that we have established that also for Erlang distributed demand an augmented wholesale price contract may not be suitable, we will examine whether a wholesale price contract with endogenous renewal probability again offers a suitable alternative. We analyze whether under the coordinating wholesale price both parties can earn a positive profit and how the total profits are divided. 
We can observe from Figure \ref{fig: division profits erlang} that under the coordinating wholesale price, the division of profits between the supplier and the OEM for different Erlang distributions with $\lambda=1$ for $\delta=0.9$ has similar characteristics as the division that was found in Figure \ref{fig:distr_profits} (with Erlang-1 being equal to the exponential distribution). More specifically, when a substantial part of demand is fixed, the distribution of profits is similar for the different values of Erlang shape parameter $n$.

\begin{figure}
    \centering
    \includegraphics[width=0.5\textwidth]{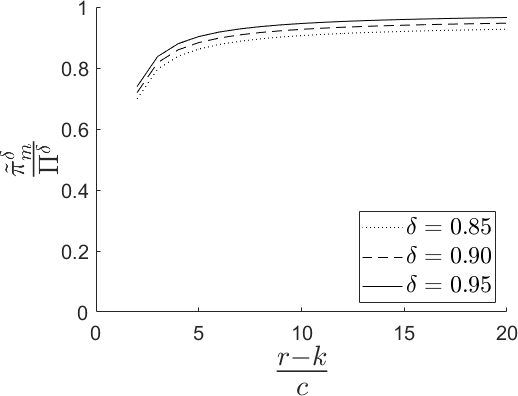}
    \caption{Fraction of the NPV captured by the OEM ($\tilde{\pi}^\delta_m/\Pi^\delta$) under the coordinating contingent wholesale price contract with Erlang-3 demand for different values of $\delta$.}
    \label{fig: division profits erlang delta}
\end{figure}

In \S\ref{sec: endogenous renewal - surplus extraction}, we found that the fraction of the profit captured by the OEM increases with $\delta$. From Figure \ref{fig: division profits erlang delta} we observe that also for $b=1$ and $A\sim Erlang(1,3)$, $\delta$ has a positive effect on the share of profit captured by the OEM. In summary, even when we extend our analysis to the case when the uncertain demand $A$ follows the $Erlang(\lambda, n)$ distribution, the structural results presented in \S4.3.2 as depicted in  Figure \ref{fig:distr_profits} continue to hold.

\section{Discussion and conclusion}\label{sec: conclusions}

Motivated by our discussions with different European OEMs in the high-tech industry, we have examined different types of supply chain contracts between a high-tech OEM who designs and manufactures multiple generations of a state-of-the-art system and the supplier of a critical component of this system.  Different from the existing literature, our model captures certain unique characteristics that are prevalent in the high-tech industry: demand consists of advance orders and a highly uncertain additional demand; components for each generation are single-sourced; capacity established by the supplier is not verifiable by the OEM; and the under-stock costs is very high, or, equivalently, the selling price is very high. Consequently, the cost of underinvestment by the supplier under a standard wholesale price contract is high and the OEM seeks for opportunities to entice the supplier to invest in more capacity.

Our work complements existing literature by examining two new supply contracts that are of practical relevance to the OEMs in the high-tech industry: (a) augmented wholesale price contracts with contingent penalty for a single generation; and (b) expanded wholesale price contracts with contingent renewal for multiple generations. 
By examining the equilibrium outcomes, we have established the following results.  First, the augmented wholesale price contract can coordinate the supply chain and it is optimal (in the sense that the OEM can capture all profit from the supply chain).  Second, the contingent renewal contract can coordinate the supply chain but the OEM cannot capture the entire supply chain profit.  Third, the augmented wholesale price contract is more efficient when the margin is very low, but the contingent renewal contract is more practical when the margin is very high. 
More specifically, in a high-tech setting, where the margin of the end product is usually large, an augmented wholesale price contract may not be enforceable in practice, while in this case the OEM can earn a large share of total profits when using a contingent renewal contract. Such a contingent renewal contract is especially attractive for the OEM when expected demand per period is large or when the supplier has a high valuation of future profits, since the supplier will be more willing to invest in capacity under these circumstances. Long-term collaborations are thus not only useful, but also essential for the functioning of high-tech supply chains.

Even though our problem setting differs from that of \textcite{tay07,taylor2007simple} in important ways, there are some similarities in the obtained results. Most importantly, both their studies and our study have shown that the gain from long-term collaboration is largest when the valuation of future profits is high. In this setting, \textcite{taylor2007simple} recommend a price-quantity contract. Our analytic result for $\frac{r-k}{c}$ large demonstrates that for $\delta$ large, even the price-only contract may perform well, which is important since such contracts are easier to adopt in the high-tech setting. We also found that the OEM can capture more profit when the amount of base demand increases. 

Additionally, we analyzed the effect of including a positive reservation profit of the supplier and of using other demand distributions on our results. When considering the supplier's reservation profit we find that even though the supplier's expected profit is equal to the reservation profit, the augmented wholesale price contract with lump-sum penalty faces the same problems as without reservation profit. Hence, the main difficulties of the augmented wholesale price contract are not countered by including a positive reservation profit. Next, we showed by means of numerical analysis that the obtained results do not only hold when demand consists of a fixed base demand and an exponential part, but extend to more general demand distributions such as Erlang distributed demand.

Even though the model presented in this study provides insights into collaborations in high-tech supply chains and shows the value of establishing long-term interactions, it has several limitations for further examination in the future. First, in our model we assume that the demand function and the cost parameters are generation-independent. However, there may be trends in demand that are not captured by this model. Therefore, it is of interest to extend our model to generation-dependent demand functions, e.g. $\bb_t= \bb_0 \alpha^t$, with $\alpha>1$ the growth factor of demand. Another interesting aspect is when production equipment (i.e. capacity investment) can be used for more than 1 product generation.

\section*{Acknowledgements}
This work is part of the research programme Complexity in high-tech manufacturing
with project number 439.16.121, which is (partly) financed by the Dutch Research Council (NWO).

\printbibliography

\input{Appendix.tex}

\end{document}

%% file: Appendix.tex
\begin{appendices}
\section{Proofs}

\paragraph{Proof of Lemma \ref{lem:helperlemma}}
\begin{proof}
Let $f(y)=(y-k-c)(b+\frac{1}{\lambda})-\frac{c}{\lambda}\ln\left(\frac{y-k}{c}\right)$. We now prove the claim that $f(y)\ge 0$ if $y\ge k+c$.  To begin, let $z = \frac{y-k}{c}-1$ so that
$f(z) = cz (b+\frac{1}{\lambda})-\frac{c}{\lambda}\ln(1+z)$.  By noting 
that $f(0) = 0$ and $f(z)$ is increasing and convex in $z \ge 0$ and by noting that $z \ge 0$ when $y\ge k+c$, we can conclude that $f(y) \ge 0$. \end{proof}



\paragraph{Proof of Proposition \ref{prop:sp_centralized}}
\begin{proof}
For any $x>\bb$ we have:
\begin{align*}
    \Pi(x) &= -cx + \expect{(r-k)\min\{\bb+\A,x\}}\\
    &= -cx+(r-k)\bb+\frac{r-k}{\lambda}\left(1-e^{-\lambda(x-\bb)}\right)
    \intertext{Taking the derivative w.r.t $x$ gives}
    \frac{d}{dx}\Pi(x)&=-c+(r-k)e^{-\lambda(x-\bb)}.
    \intertext{Thus, using that $r-k>2c$, we find that $\frac{d}{dx}\Pi(x)=0$ iff:}
    x&=\bb+\frac{1}{\lambda}\ln\left(\frac{r-k}{c}\right).
    \intertext{Hence, since $(r-k)/c>1$, this implies that the $x$ that maximizes $\Pi(x)$ must satisfy $x^*=\bb+\frac{1}{\lambda}\ln\left(\frac{r-k}{c}\right)$ Also, corresponding supply chain profit is:}
    \Pi^*&=-c\left(\bb+\frac{1}{\lambda}\ln\left(\frac{r-k}{c}\right)\right)+(r-k)\bb + \frac{r-k}{\lambda}\left(1-\frac{c}{r-k}\right)\\
    &=(r-k-c)\left(\bb+\frac{1}{\lambda}\right)-\frac{c}{\lambda}\ln\left(\frac{r-k}{c}\right).
\end{align*}
\textcolor{black}{
Now, to prove that $\Pi^*$ is positive, note that  $\Pi^*=(r-k-c)(b+\frac{1}{\lambda})-\frac{c}{\lambda}\ln\left(\frac{r-k}{c}\right)=f(r)$, with $f(\cdot)$ as in the proof of Lemma~\ref{lem:helperlemma}. From the lemma and the assumption $r>k+c$, it then follows that $\Pi^*>0$. 
}
\end{proof}

\paragraph{Proof of Proposition \ref{prop:ws}}
\begin{proof}
Taking into account the supplier's capacity decision, the OEM's profit is:
\begin{align*}
    \tilde{\pi}_m(w) & = (r-w)\bb+\frac{r-w}{\lambda}\left(1-\frac{c}{w-k}\right)
    \intertext{of which the derivative w.r.t $w$ equals}
    \frac{d}{dw}\tilde{\pi}_m(w) &= -\left(\bb+\frac{1}{\lambda}\right)+\frac{1}{\lambda}\frac{(r-k)c}{(w-k)^2}
    \intertext{Clearly, $\frac{d^2}{dw^2}\tilde{\pi}(w)<0$ for $w>k$ and thus $\tilde{\pi}_m(w)$ is maximized by setting the first derivative equal to $0$, yielding}
    \tilde{w}&= k+\sqrt{\frac{(r-k)c}{\bb\lambda+1}}.
    \intertext{Since $r-k-c>c \rightarrow r-k>2c$, and since $\frac{1}{\lambda}>b\rightarrow b\lambda+1 < 2$, we find $\tilde{w}=k+\sqrt{\frac{(r-k)c}{b\lambda+1}}>k+\sqrt{2c^2/2}=k+c$, thus the supplier participation constraint is satisfied. The resulting capacity developed by the supplier follows by substituting this in \eqref{eq:sp_dc_x}}
    \tilde{x}&=\bb+\frac{1}{\lambda} \ln\left(\sqrt{\frac{r-k}{(\bb\lambda+1)c}}\right).
\end{align*}
\end{proof}

\paragraph{Proof of Proposition \ref{prop:ws_p}}
\begin{proof}
\textcolor{black}{We first prove that the manufacturers and suppliers profit are strictly positive. Note that the manufacturer's profit is obtained by substituting $\tilde{w}$ into the manufacturer's profit function $\tilde{\pi}_m(w)$. Since $r-k-c>c \rightarrow r-k>2c$, and since $\frac{1}{\lambda}>b\rightarrow b\lambda+1 < 2$, we find $\tilde{w}=k+\sqrt{\frac{(r-k)c}{b\lambda+1}}>k+\sqrt{2c^2/2}=k+c$. Since $r-k\le r$, $c<r$, and $b\lambda+1\ge 1$, we also have $\tilde{w}<r$. Thus $\tilde{\pi}_m=\tilde{\pi}_m(\tilde{w})=(r-\tilde{w})(b+\frac{1}{\lambda}(1-\frac{c}{\tilde{w}-k})>0$. The suppliers profit is obtained by substituting $\tilde{w}$ in $\tilde{\pi}_s(w)$. Note that $\tilde{\pi}_s(\tilde{w})=f(\tilde{w})$, with $f(\cdot)$ as in the proof of Lemma~\ref{lem:helperlemma}. That $\tilde{\pi}_s>0$ then follows from Lemma~\ref{lem:helperlemma} and because $\tilde{w}>k+c$, as was shown above. }

To prove that the total profit is lower than in the centralized case, we substitute the supplier's capacity and OEM's wholesale price decisions provided in Proposition \ref{prop:ws} in the supplier's and OEM's profit functions given in Equations \eqref{pisd} and \eqref{pimd}. This confirms our expressions in the Proposition:
\begin{align*}
    \tilde{\pi}_s &= \left(\sqrt{\frac{(r-k)c}{\bb\lambda+1}}-c\right)\bb + \frac{\sqrt{\frac{(r-k)c}{\bb\lambda+1}}}{\lambda}-\frac{c}{\lambda}-\frac{c}{\lambda}\ln\left(\sqrt{\frac{r-k}{(\bb\lambda+1)c}} \right), \mbox{ and}\\
    \tilde{\pi}_m &=\left(r-k-\sqrt{\frac{(r-k)c}{\bb\lambda+1}}\right) \left(\bb + \frac{1}{\lambda}\left(1-\sqrt{\frac{(\bb\lambda+1)c}{r-k}}\right)\right).
\end{align*}
Supply chain profit thus equals
\begin{align*}
    \tilde{\pi}_s+\tilde{\pi}_m = &\left(\sqrt{\frac{(r-k)c}{\bb\lambda+1}}-c\right)\bb + \frac{\sqrt{\frac{(r-k)c}{\bb\lambda+1}}}{\lambda}-\frac{c}{\lambda}-\frac{c}{\lambda}\ln\left(\sqrt{\frac{r-k}{(\bb\lambda+1)c}} \right)\\
    &+ \left(r-k-\sqrt{\frac{(r-k)c}{\bb\lambda+1}}\right) \left(\bb + \frac{1}{\lambda}\left(1-\sqrt{\frac{(\bb\lambda+1)c}{r-k}}\right)\right)\\
    =& (r-k-c)\left(\bb+\frac{1}{\lambda}\right) -\frac{c}{\lambda} \ln\left(\sqrt{\frac{r-k}{(\bb\lambda+1)c}}\right)-\frac{r-k-\sqrt{\frac{(r-k)c}{\bb\lambda+1}}}{\lambda}\sqrt{\frac{(\bb\lambda+1)c}{r-k}}\\
    =& (r-k-c)\left(\bb+\frac{1}{\lambda}\right) -\frac{c}{\lambda} \ln\left(\sqrt{\frac{r-k}{(\bb\lambda+1)c}}\right)- \frac{1}{\lambda}\sqrt{(r-k)(\bb\lambda+1)c}+\frac{c}{\lambda} 
\end{align*}
Since
\begin{align*}
    \frac{c}{\lambda} \ln\left(\sqrt{\frac{r-k}{(\bb\lambda+1)c}}\right) + \frac{1}{\lambda}\sqrt{(r-k)(\bb\lambda+1)c}-\frac{c}{\lambda} > \frac{c}{\lambda}\ln\left(\frac{r-k}{c}\right)
    \intertext{or equivalently}
    \sqrt{\frac{(r-k)(\bb\lambda+1)}{c}} > 1 + \ln\left(\sqrt{\frac{(r-k)(\bb\lambda+1)}{c}}\right)
\end{align*}
it follows that $ \tilde{\pi}_s+\tilde{\pi}_m < \Pi^*$.
\end{proof}

\paragraph{Proof of Proposition \ref{prop:penaltycontract}}


\begin{proof}
Suppose that the supplier participates in the penalty contract. Then his profit equals:
\begin{align*}
    \hat{\pi}_s(x,w,\rho)&= -cx+ \expect{(w-k) \min(\bb+\A, x) - \rho 1_{\{\bb+\A > x\}}}\\
    &= (w-c-k)\bb-c(x-\bb)+\frac{w-k}{\lambda}\left(1-e^{-\lambda(x-\bb)}\right)-\rho e^{-\lambda(x-\bb)}.
    \intertext{By taking the derivative with respect to $x$, and setting it equal to zero, we find the capacity level at which the supplier maximizes the profit of participating in the contract:}
    \hat{x}(w,\rho)&=\bb+\frac{1}{\lambda}\ln\left(\frac{w-k+\rho\lambda}{c}\right).
    \intertext{Now note that by Proposition~\ref{prop:sp_centralized}, the supply chain profit is maximized if $x=\bb+\frac{1}{\lambda}\ln(\frac{r-k}{c})$. In view of the above, a necessary condition to achieve this is that $\frac{w-k+\rho\lambda}{c}=\frac{r-k}{c} \rightarrow \rho\lambda+w=r$. To arrive at a coordinating contract in which the OEM captures all the profit, we must in addition set $\rho$ and $w$ such that the profit of participating for the supplier equals $0$. The supplier's profit equals:
    }
    \hat{\pi}_s(\hat{x}(w,\rho),w,\rho)&= (w-c-k)\bb+\frac{w-k-c}{\lambda}-\frac{c}{\lambda}\ln\left(\frac{w-k+\rho\lambda}{c}\right)
    \intertext{we impose $\rho\lambda+w=r$ by substituting $w=r-\rho\lambda$, which yields}
    \hat{\pi}_s(\rho)&= (r-\rho\lambda-c-k)\bb+\frac{r-\rho\lambda-k-c}{\lambda}-\frac{c}{\lambda}\ln\left(\frac{r-k}{c}\right)
    \intertext{We in addition impose $\hat{\pi}_s(\rho)$=0, which holds for:}
    \hat{\rho} &= \frac{1}{\lambda}\left(r-k-c-\frac{c}{\bb\lambda+1}\ln\left(\frac{r-k}{c}\right)\right).
    \intertext{Note that $\hat{\rho}>0$ whenever $r-k-c>\frac{c}{\bb\lambda+1}\ln\left(\frac{r-k}{c}\right)$, or equivalently $\frac{r-k}{c}>1+\frac{1}{\bb\lambda+1}\ln\left(\frac{r-k}{c}\right)$. For $\frac{r-k}{c}=1$ both sides would equal 1. For $\frac{r-k}{c}>1$ the left-hand side increases linearly in $\frac{r-k}{c}$, while the right-hand side increases logarithmically. Also, since $\bb\lambda>0$, the fraction $\frac{1}{\bb\lambda+1}<1$. Therefore, using our earlier assumptions on $r, \bb$ and $\frac{1}{\lambda}$ it holds that $\frac{r-k}{c}>1+\frac{1}{\bb\lambda+1}\ln\left(\frac{r-k}{c}\right)$ and it thus follows that $\hat{\rho}>0$. To satisfy $r=w+\rho\lambda$, the OEM sets}
    \hat{w}&=k+c+\frac{c}{\bb\lambda+1}\ln\left(\frac{r-k}{c}\right).
    \intertext{By construction $\hat{\pi}_s(\hat{x}(\hat{w},\hat{\rho}),\hat{w},\hat{\rho})=0$ and it follows that}
    \hat{\pi}_m(\hat{w},\hat{\rho}) &= \left(r-k-c-\frac{c}{\bb\lambda+1}\ln\left(\frac{r-k}{c}\right)\right)\left(\bb+\frac{1}{\lambda}\right)\\
    &= \Pi^*.
\end{align*}
\end{proof}

\paragraph*{Proof of Proposition \ref{prop:sup_opt}}
\begin{proof}
Let $y=e^{-\lambda (x-\bb)}$. Then we can rewrite Equation \eqref{eq:mp_s_p} as 
\begin{equation*}
\tilde{\pi}_s^\delta(y)=\frac{1}{1-\delta(1-y)}\left(\bb(w-k-c) + \frac{c}{\lambda}\ln(y)+\frac{w-k}{\lambda}(1-y)\right)
\end{equation*}
The first derivative with respect to $y$ is given by
\begin{equation*}
\frac{d}{dy}\tilde{\pi}_s^\delta(y)=\frac{(1-\delta(1-y))\left(\frac{c}{\lambda}\frac{1}{y}-\frac{w-k}{\lambda}\right)-\left(\bb(w-k-c)+\frac{c}{\lambda}\ln(y)+\frac{w-k}{\lambda}(1-y)\right)\delta}{(1-\delta(1-y))^2}
\end{equation*}
Equating the derivative to 0 yields
\begin{align*}
(1-\delta(1-y))\left(\frac{c}{\lambda}\frac{1}{y}-\frac{w-k}{\lambda}\right)&=\left(\bb(w-k-c)+\frac{c}{\lambda}\ln(y)+\frac{w-k}{\lambda}(1-y)\right)\delta
\intertext{which can be rewritten as}
\left(w-k-\delta c + \delta \bb(w-k-c)\lambda\right)y + \delta c y\: \ln(y) &= (1-\delta)c
\intertext{which is of the form}
a_1 y + a_2 y \ln(y) &= a_3
\intertext{with $a_1,a_2,a_3>0$ since $w\ge c+k$, $b\ge 0$, $\lambda>0$ and $0<\delta<1$. Therefore, we can rewrite this as}
e^{a_1 y + a_2 y \ln(y)} &= e^{a_3}\\
e^{a_2 y \ln(y)} &= e^{a_3-a_1 y}\\
e^{\ln(y)} &= e^{\frac{a_3-a_1 y}{a_2 y}}\\
y &= e^{\frac{a_3}{a_2 y}-\frac{a_1}{a_2}}\\
y e^{\frac{a_1}{a_2}} &= e^{\frac{a_3}{a_2 y}}\\
\frac{a_3}{a_2}e^{\frac{a_1}{a_2}} &= \frac{a_3}{a_2 y} e ^{\frac{a_3}{a_2 y}}\\
W\left(\frac{a_3}{a_2}e^\frac{a_1}{a_2}\right) &= \frac{a_3}{a_2 y}
\intertext{to reach solution}
y &= \frac{a_3}{a_2 W\left(\frac{a_3}{a_2}e^\frac{a_1}{a_2}\right)}
\intertext{where $W$ is the Lambert's W function. Substituting $a_1=w-k-\delta c + \delta \bb(w-k-c)\lambda$, $a_2=\delta c$ and $a_3 = (1-\delta)c$ and simplifying the expression yields}
y &=\frac{1-\delta}{\delta W\left(\frac{1-\delta}{\delta}e^{\frac{w-k}{\delta c}-1+\frac{\lambda \bb(w-k-c)}{c}}\right)}
\end{align*}
Since $y=e^{-\lambda (x-\bb)}$, we obtain
\begin{equation*}
\tilde{x}(w,\delta)=\bb+\frac{1}{\lambda}\ln\left(\frac{\delta}{1-\delta}W\left(\frac{1-\delta}{\delta}e^{\frac{w-k}{\delta c}-1+\frac{\lambda \bb(w-k-c)}{c}}\right)\right)
\end{equation*}
For the supplier participation constraint, note that Lemma~\ref{lem:helperlemma} proved that setting capacity $\tilde{x}(w)$ when offered a wholesale price $w>k+c$ yields positive expected profits $\tilde{\pi}_s(\tilde{x}(w),w)$ per generation for the supplier whenever $w>k+c$. Hence, for any $w>k+c$, we must have that $\tilde{\pi}^\delta_s(\tilde{x}(w))=\tilde{\pi}_s(\tilde{x}(w),w)/(1-\delta R(\tilde{x}(w)))$. This shows that for any $w>k+c$, there exists a capacity level that yields positive profits for the supplier, hence the supplier will participate in any contingent wholesale contract with $w>k+c$. 
\end{proof}

The main text claims the following, and we present a formal proof here for completeness.
\begin{proposition}
For given $w$, $\tilde{x}(w,\delta)$ is increasing in $\delta\in (0,1)$. \label{prop:x_incr_delta}
\end{proposition}

\paragraph*{Proof of Proposition \ref{prop:x_incr_delta}}
\begin{proof}
Let $Z=\frac{\lambda \bb(w-k-c)}{c}$. Then 
\begin{equation*}
\tilde{x}(w,\delta)=\bb+\frac{1}{\lambda}\ln\left(\frac{\delta}{1-\delta}W\left(\frac{1-\delta}{\delta}e^{\frac{w-k}{\delta c}-1+Z}\right)\right)
\end{equation*} 
For $\delta\in(0,1)$ it is readily seen that both $\frac{1-\delta}{\delta}$ and $e^{\frac{w-k}{\delta c}-1+Z}$ are strictly decreasing in $\delta$, thus so is $\frac{1-\delta}{\delta} e^{\frac{w-k}{\delta c}-1+Z}$. Combined with the fact that $W(x)$ is increasing in $x\geq0$ this shows that $W\left(\frac{1-\delta}{\delta}e^{\frac{w-k}{\delta c}-1+Z}\right)$ is decreasing in $\delta$. 

Considering $\frac{\delta}{1-\delta}W\left(\frac{1-\delta}{\delta}e^{\frac{w-k}{\delta c}-1+Z}\right)$, this thus consists of a term $f(\delta)=\frac{\delta}{1-\delta}$ that is strictly increasing in $\delta$ and the term $g(\delta)=W\left(\frac{1-\delta}{\delta}e^{\frac{w-k}{\delta c}-1+Z}\right)$ that is strictly decreasing in $\delta$, thus $f'(\delta)g'(\delta)<0$. Therefore, the entire term is increasing in $\delta$ if and only if
\begin{equation*}
 (f(\delta)g(\delta))'=f'(\delta)g(\delta)+f(\delta)g'(\delta)>0\leftrightarrow \frac{g(\delta)}{g'(\delta)}< - \frac{f(\delta)}{f'(\delta)}
\end{equation*}
Substituting 
\begin{align*}
f'(\delta)&= \frac{-1}{(1-\delta)^2}\\
g'(\delta) &= -\frac{W\left(\frac{1-\delta}{\delta}e^{\frac{w-k}{\delta c}-1+Z}\right)}{W\left(\frac{1-\delta}{\delta}e^{\frac{w-k}{\delta c}-1+Z}\right)+1}\left(\frac{w-k}{\delta^2c}+\frac{1}{\delta(1-\delta)}\right)
\end{align*}
yields the equivalent inequality
\begin{equation*}
\frac{W\left(\frac{1-\delta}{\delta}e^{\frac{w-k}{\delta c}-1+Z}\right)+1}{\frac{w-k}{\delta^2c}+\frac{1}{\delta(1-\delta)}} > \delta(1-\delta)
\end{equation*}
which can be rewritten as
\begin{equation*}
\frac{w-k}{c}+Z> \ln\left(\frac{w-k}{c}\right)+1
\end{equation*}
Since $Z\geq0$, this holds for all $w>c+k$. 
Therefore, $\frac{\delta}{1-\delta}W\left(\frac{1-\delta}{\delta}e^{\frac{w-k}{\delta c}-1+Z}\right)$ is increasing in $\delta$ and so is $\ln\left(\frac{\delta}{1-\delta}W\left(\frac{1-\delta}{\delta}e^{\frac{w-k}{\delta c}-1+Z}\right)\right)$. 
Thus we can conclude that $\tilde{x}(w,\delta)$ is increasing in $\delta\in(0,1)$.
\end{proof}

\paragraph*{Proof of Proposition \ref{prop:sl}}
\begin{proof}
Under the single-epoch wholesale price contract, the supplier's capacity decision equals
\begin{equation*}
    \tilde{x}(w)=\bb+\frac{1}{\lambda}\ln\left(\frac{w-k}{c}\right)
\end{equation*}
Under the long-term contract with performance contingency, the supplier's capacity decision equals
\begin{equation*}
    \tilde{x}(w,\delta)=\bb+\frac{1}{\lambda}\ln\left(\frac{\delta}{1-\delta}W\left(\frac{1-\delta}{\delta}e^{\frac{w-k}{\delta c}-1+\frac{\lambda \bb(w-k-c)}{c}}\right)\right).
\end{equation*}
It follows that $\tilde{x}(w,\delta)>\tilde{x}(w)$ if and only if
\begin{equation*}
    \frac{\delta}{1-\delta}W\left(\frac{1-\delta}{\delta}e^{\frac{w-k}{\delta c}-1+\frac{\lambda \bb(w-k-c)}{c}}\right) > \frac{w-k}{c}
\end{equation*}
This relationship holds if and only if
\begingroup
\allowdisplaybreaks
\begin{alignat*}{3}
&& W\left(\frac{1-\delta}{\delta}e^{\frac{w-k}{\delta c}-1+\frac{\lambda \bb(w-k-c)}{c}}\right) & > \frac{w-k}{c}\frac{1-\delta}{\delta}\\
\Leftrightarrow\hspace{1cm} &&\frac{1-\delta}{\delta}e^{\frac{w-k}{\delta c}-1+\frac{\lambda \bb(w-k-c)}{c}} & > W^{-1}\left(\frac{w-k}{c}\frac{1-\delta}{\delta}\right)\\
\Leftrightarrow\hspace{1cm} &&\frac{1-\delta}{\delta}e^{\frac{w-k}{\delta c}-1+\frac{\lambda \bb(w-k-c)}{c}} & > \frac{w-k}{c}\frac{1-\delta}{\delta}e^{\frac{w-k}{c}\frac{1-\delta}{\delta}}\\
\Leftrightarrow\hspace{1cm} && e^{\frac{w-k}{\delta c}-1+\frac{\lambda \bb(w-k-c)}{c}} & > \frac{w-k}{c} e^{\frac{w-k}{c}\frac{1-\delta}{\delta}}\\
\Leftrightarrow\hspace{1cm} && \frac{w-k}{\delta c}-1+\frac{\lambda \bb(w-k-c)}{c} & > \ln\left(\frac{w-k}{c}\right)+\frac{w-k}{c}\frac{1-\delta}{\delta}\\
\Leftrightarrow\hspace{1cm} &&  \frac{\lambda \bb(w-k-c)}{c}-1 & > \ln\left(\frac{w-k}{c}\right)-\frac{w-k}{c}\\
\Leftrightarrow\hspace{1cm} &&  \frac{\lambda \bb(w-k-c)}{c}+\frac{w-k}{c}-1 & > \ln\left(\frac{w-k}{c}\right)
\end{alignat*}
\endgroup
Here, in the first equivalence, we use that $W(\cdot)$ is strictly increasing, for the second equivalence, we used the definition of $W(\cdot)$, while the third equivalence is obtained by taking the logarithm on both sides. The remaining equivalences are simple manipulations. This holds for all $w-k>c$ and $\bb\geq 0$. 
\end{proof}

\paragraph{Proof of Proposition \ref{prop:mp_dc_c_w}}
\begin{proof}
The supply chain is coordinated when the OEM sets $w$ such that $\tilde{x}(w,\delta)=x^*$. This is the case when 
\begin{equation*}
    \frac{\delta}{1-\delta}W\left(\frac{1-\delta}{\delta}e^{\frac{w-k}{\delta c}-1+\frac{\lambda \bb(w-k-c)}{c}}\right)=\frac{r-k}{c}
\end{equation*}
This equality holds if and only if
\begingroup
\allowdisplaybreaks
\begin{alignat*}{3}
&& \frac{1-\delta}{\delta}e^{\frac{w-k}{\delta c}-1+\frac{\lambda \bb(w-k-c)}{c}} & =  W^{-1}\left(\frac{1-\delta}{\delta}\frac{r-k}{c}\right)\\
\Leftrightarrow\hspace{1cm} &&\frac{1-\delta}{\delta}e^{\frac{w-k}{\delta c}-1+\frac{\lambda \bb(w-k-c)}{c}} & = \frac{1-\delta}{\delta}\frac{r-k}{c}e^{\frac{1-\delta}{\delta}\frac{r-k}{c}}\\
\Leftrightarrow\hspace{1cm} && e^{\frac{w-k}{\delta c}-1+\frac{\lambda \bb(w-k-c)}{c}} & = \frac{r-k}{c}e^{\frac{r-k}{\delta c}-\frac{r-k}{c}}\\
\Leftrightarrow\hspace{1cm} && \frac{w-k}{\delta c}-1+\frac{\lambda \bb(w-k-c)}{c} & = \ln\left(\frac{r-k}{c}\right)+ \frac{r-k}{\delta c}-\frac{r-k}{c}\\
\Leftrightarrow\hspace{1cm} && w-k-\delta c+ \delta\lambda \bb(w-k-c) & = \delta c\: \ln\left(\frac{r-k}{c}\right)+r-k-\delta(r-k)\\
\Leftrightarrow\hspace{1cm} &&  (w-k)(1+\delta \bb\lambda) -\delta c (1+\bb\lambda) & = \delta c\: \ln\left(\frac{r-k}{c}\right)+(1-\delta)(r-k)\\
\Leftrightarrow\hspace{1cm} &&  w & = k+\frac{\delta c \left(1+\bb\lambda+\ln\left(\frac{r-k}{c}\right)\right)+(1-\delta)(r-k)}{1+\delta \bb\lambda}
\end{alignat*}
\endgroup
Thus $w^\delta = k+\frac{\delta c \left(1+\bb\lambda+\ln\left(\frac{r-k}{c}\right)\right)+(1-\delta)(r-k)}{1+\delta \bb\lambda}$. The resulting profits for the supplier and the OEM are:
\begin{align*}
    \tilde{\pi}_s=&(w^\delta-c-k)\bb+\frac{w^\delta-k}{\lambda}\left(1-\frac{c}{r-k}\right)-\frac{c}{\lambda}\ln\left(\frac{r-k}{c}\right)\\
    =& \frac{(1-\delta)(r-k-c)+\delta c\:\ln\left(\frac{r-k}{c}\right)}{1+\delta \bb\lambda}\bb+\frac{1}{\lambda}\frac{\delta c \left(1+\bb\lambda+\ln\left(\frac{r-k}{c}\right)\right)+(1-\delta)(r-k)}{1+\delta \bb\lambda}\left(1-\frac{c}{r-k}\right)-\frac{c}{\lambda}\ln\left(\frac{r-k}{c}\right)\\
    \intertext{and}
    \tilde{\pi}_{m}=&(r-w^\delta)\bb+\frac{r-w^\delta}{\lambda}\left(1-\frac{c}{r-k}\right)\\
    =& \frac{\delta \bb\lambda(r-k-c)+\delta(r-k-c)-\delta c\:\ln\left(\frac{r-k}{c}\right)}{1+\delta \bb\lambda}\bb +\frac{1}{\lambda}\frac{\delta \bb\lambda(r-k-c)+\delta(r-k-c)-\delta c\:\ln\left(\frac{r-k}{c}\right)}{1+\delta \bb\lambda}\left(1-\frac{c}{r-k}\right)
\end{align*}
Since $w-k>c$, the supplier's participation constraint is satisfied and $\tilde{\pi}_s>0$. 

We still need to prove $w^\delta<r$, which is equivalent to
\begin{equation*}
   \frac{\delta c \left(1+\bb\lambda+\ln\left(\frac{r-k}{c}\right)\right)+(1-\delta)(r-k)}{1+\delta \bb\lambda}< r-k
\end{equation*}
This equality holds if and only if
\begingroup
\allowdisplaybreaks
\begin{alignat*}{3}
&& \delta c \left(1+\bb\lambda+\ln\left(\frac{r-k}{c}\right)\right)+(1-\delta)(r-k) & < (1+\delta \bb\lambda)(r-k)\\
\Leftrightarrow\hspace{1cm} &&\delta c \left(1+\bb\lambda+\ln\left(\frac{r-k}{c}\right)\right)-\delta(r-k) & < \delta \bb\lambda(r-k)\\
\Leftrightarrow\hspace{1cm} && 1+\bb\lambda+\ln\left(\frac{r-k}{c}\right) - \frac{r-k}{c} & < \bb\lambda\frac{r-k}{c}\\
\Leftrightarrow\hspace{1cm} && 1+\bb\lambda+\ln\left(\frac{r-k}{c}\right) & < (\bb\lambda+1)\frac{r-k}{c}
\end{alignat*}
\endgroup
For $r-k>c$, since $\bb\lambda\geq 0$ the left-hand side increases logarithmically in $\frac{r-k}{c}$, while the right-hand side increases linearly in $\frac{r-k}{c}$. Therefore, we conclude that the inequality holds and the coordinating wholesale price when considering multiple product generations is lower than $r$.
\end{proof}



\paragraph{Proof of Proposition \ref{prop:NPVcapture}}
\begin{proof}
If the fraction of the NPV of the total supply chain profit captured by the OEM is smaller than 1, this means that   $\tilde{\pi}_{m}^{\delta} < \Pi^\delta$, or equivalently $\Pi^\delta-\tilde{\pi}_{m}^{\delta}>0$.
Inserting $\tilde{\pi}_{m}^{\delta}
=  \frac{1}{1-\delta}\left(r-w^\delta\right)\left(\bb+\frac{1}{\lambda}(1-\frac{c}{r-k})\right)$ with $ w^\delta = k+\frac{\delta c \left(1+\bb\lambda+\ln\left(\frac{r-k}{c}\right)\right)+(1-\delta)(r-k)}{1+\delta \bb\lambda}$ and $\Pi^{\delta} =\frac{1}{1-\delta}\left((r-k-c)\left(\bb+\frac{1}{\lambda}\right) - \frac{c}{\lambda}\ln\left(\frac{r-k}{c}\right)\right)$ and simplifying the resulting expression yields the condition:
\begin{equation*}
    \frac{(c\delta+(1-\delta)(r-k))\left(\left(b+\frac{1}{\lambda}\right)(r-c-k)-c\: ln\left(\frac{r-k}{c}\right)\right)}{(1+\delta b\lambda)(r-k)}>0.
\end{equation*}
Since $r>k$ and $\delta<1$, this holds when $\left(b+\frac{1}{\lambda}\right)(r-c-k)>c \: ln\left(\frac{r-k}{c}\right)$, or equivalently $\frac{r-k}{c}>1+\frac{1}{b\lambda+1}ln\left(\frac{r-k}{c}\right)$. This holds for all $r-k-c>0$, by the same reasoning as in the proof of Proposition \ref{prop:mp_dc_c_w}.
\end{proof}

\paragraph{Proof of Proposition \ref{prop:distr_profits}}
\begin{proof}
The OEM's profit is denoted as in Proposition \ref{prop:mp_dc_c_w} by
\begin{align*}
    \tilde{\pi}_{m}=& \frac{\delta \bb\lambda(r-k-c)+\delta(r-k-c)-\delta c\:\ln\left(\frac{r-k}{c}\right)}{1+\delta \bb\lambda}\bb+\frac{1}{\lambda}\frac{\delta \bb\lambda(r-k-c)+\delta(r-k-c)-\delta c\:\ln\left(\frac{r-k}{c}\right)}{1+\delta \bb\lambda}\left(1-\frac{c}{r-k}\right)
\end{align*}
This means that the fraction of supply chain profit captured by the OEM equals
\begin{align*}
    \frac{ \tilde{\pi}_{m}}{\Pi^*}&=\frac{\frac{\delta \bb\lambda(r-k-c)+\delta(r-k-c)-\delta c\:\ln\left(\frac{r-k}{c}\right)}{1+\delta \bb\lambda}\bb +\frac{1}{\lambda}\frac{\delta \bb\lambda(r-k-c)+\delta(r-k-c)-\delta c\:\ln\left(\frac{r-k}{c}\right)}{1+\delta \bb\lambda}\left(1-\frac{c}{r-k}\right)}{(r-k-c)\left(\bb+\frac{1}{\lambda}\right)-\frac{c}{\lambda}\ln\left(\frac{r-k}{c}\right)}
    \intertext{or equivalently, by dividing both the numerator and denominator by $r-k$}
    &=\frac{\frac{\delta \bb\lambda\frac{r-k-c}{r-k}+\delta\frac{r-k-c}{r-k}-\delta \frac{c}{r-k}\:\ln\left(\frac{r-k}{c}\right)}{1+\delta \bb\lambda}\bb +\frac{1}{\lambda}\frac{\delta \bb\lambda\frac{r-k-c}{r-k}+\delta\frac{r-k-c}{r-k}-\delta \frac{c}{r-k}\:\ln\left(\frac{r-k}{c}\right)}{1+\delta \bb\lambda}\frac{r-k-c}{r-k}}{\frac{r-k-c}{r-k}\left(\bb+\frac{1}{\lambda}\right)-\frac{c}{r-k}\frac{1}{\lambda}\ln\left(\frac{r-k}{c}\right)}
\end{align*}
If $\frac{r-k}{c}\rightarrow \infty$:
\begin{equation*}
    \frac{c}{r-k}\rightarrow 0, 
    \frac{r-k-c}{r-k}\rightarrow1, \textrm{ and }
    \frac{c}{r-k}\ln\left(\frac{r-k}{c}\right)\rightarrow0
\end{equation*}
Therefore,
\begin{align*}
     \frac{ \tilde{\pi}_{m}}{\Pi^*}&\rightarrow \frac{\frac{\delta \bb \lambda +\delta}{1+\delta \bb\lambda}\bb+\frac{1}{\lambda}\frac{\delta \bb \lambda + \delta}{1+\delta \bb \lambda}}{\bb+\frac{1}{\lambda}}=\frac{\delta \bb \lambda +\delta}{\delta \bb\lambda + 1}
\end{align*}
\end{proof}

\paragraph{Proof of Proposition \ref{prop:distr_profits_small_y}}
\begin{proof}
For $\delta\rightarrow1$ we obtain $w^\delta\approx k+\frac{c\left(1+b\lambda+ln\left(\frac{r-k}{c}\right)\right)}{1+b\lambda}$. As $\frac{r-k}{c}\rightarrow 1^+$, $\:ln\left(\frac{r-k}{c}\right)\rightarrow 0^+$. Therefore, $w^\delta\approx k+\frac{c\left(1+b\lambda\right)}{1+b\lambda}=k+c$. Since $\frac{r-k}{c}\rightarrow 1^+$ means that $r-(k+c)\rightarrow 0^+$ this means that $w^\delta\approx r$. 
Furthermore, in Proposition \ref{prop:penaltycontract} it is defined that $\hat{\rho} = \frac{1}{\lambda}\left(r-k-c-\frac{c}{\bb\lambda+1}\ln\left(\frac{r-k}{c}\right) \right)$. Again using $\frac{r-k}{c}\rightarrow 1^+$ and thus $\:ln\left(\frac{r-k}{c}\right)\approx 0$ we obtain $\hat{\rho} \approx \frac{1}{\lambda}(0-\frac{c}{b\lambda+1}\cdot 0)=0$.
\end{proof}

\paragraph{Proof of Proposition \ref{prop: duration collaboration coordination}}
\begin{proof}
From $Y\sim Geom(1-R(x))$ it follows that $\mathbb{E}[Y(x)]=\frac{1}{1-R(x)}$. Since $R(x)=1-e^{-\lambda(x-b)}$, we have $\mathbb{E}[Y(x)]=\frac{1}{e^{-\lambda(x-b)}}=e^{\lambda(x-b)}$. Under the coordinating contract with capacity $\tilde{x}(w^\delta)=x^*=b+\frac{1}{\lambda}\ln\left(\frac{r-k}{c}\right)$, it follows that $\mathbb{E}[Y(x)]=\frac{r-k}{c}$.
\end{proof}

\paragraph{Proof of Lemma \ref{lem: reservation profit}}
\begin{proof}
We want to determine the minimum wholesale price for which  $\tilde{\pi}_s(w)=(w-c-k)\left(b+\frac{1}{\lambda}\right)-\frac{c}{\lambda}\ln\left(\frac{w-k}{c}\right)\geq\reservationprofit$.
This holds if and only if 
\begin{align*}
    && \left(\frac{w-k}{c}-1\right)(b\lambda+1)-Z\frac{\lambda}{c} &\geq \ln\left(\frac{w-k}{c}\right)\\
    \Leftrightarrow\hspace{1cm} && e^{\frac{w-k}{c}(b\lambda+1)-(b\lambda+1+Z\frac{\lambda}{c})} &\geq \frac{w-k}{c}\\
    \Leftrightarrow\hspace{1cm} && e^{-(b\lambda+1+Z\frac{\lambda}{c})} &\geq \frac{w-k}{c}e^{-(b\lambda+1)\frac{w-k}{c}}\\
    \Leftrightarrow\hspace{1cm} && -(b\lambda+1) e^{-(b\lambda+1+Z\frac{\lambda}{c})} &\geq -(b\lambda+1)\frac{w-k}{c}e^{-(b\lambda+1)\frac{w-k}{c}}\\
    \Leftrightarrow\hspace{1cm} && W\left( -(b\lambda+1) e^{-(b\lambda+1+Z\frac{\lambda}{c})} \right) &\geq -(b\lambda+1)\frac{w-k}{c}\\
    \Leftrightarrow\hspace{1cm} && \frac{w-k}{c} &\geq -\frac{1}{b\lambda+1}W\left( -(b\lambda+1) e^{-(b\lambda+1+Z\frac{\lambda}{c})} \right)\\
    \Leftrightarrow\hspace{1cm} && w &\geq k -\frac{c}{b\lambda+1}W\left( -(b\lambda+1) e^{-(b\lambda+1+Z\frac{\lambda}{c})} \right)\\
\end{align*}
So we conclude that $\tilde{\pi}_s(w)=(w-c-k)\left(b+\frac{1}{\lambda}\right)-\frac{c}{\lambda}\ln\left(\frac{w-k}{c}\right)\geq\reservationprofit$
if and only if $w\geq k-\frac{c}{b\lambda+1}W\left(-(b\lambda+1)e^{-(b\lambda+1+\reservationprofit\frac{\lambda}{c}}\right)$.
\end{proof}

\paragraph{Proof of Proposition \ref{prop:penaltycontract reservation profit}}
\begin{proof}
Similar to the proof of Proposition \ref{prop:penaltycontract}, we have as necessary condition to achieve coordination that $\rho\lambda+w=r$. To arrive at a coordinating contract in which the supplier's profit is equal to the reservation profit $\reservationprofit$ and the OEM captures the remainder of the profit, we must in addition set $\rho$ and $w$ such that the profit of participating for the supplier equals $\reservationprofit$.

The supplier's profit equals:
\begin{align*}
    \hat{\pi}_s(\hat{x}(w,\rho),w,\rho)&= (w-c-k)\bb+\frac{w-k-c}{\lambda}-\frac{c}{\lambda}\ln\left(\frac{w-k+\rho\lambda}{c}\right)
    \intertext{we impose $\rho\lambda+w=r$ by substituting $w=r-\rho\lambda$, which yields}
    \hat{\pi}_s(\rho)&= (r-\rho\lambda-c-k)\bb+\frac{r-\rho\lambda-k-c}{\lambda}-\frac{c}{\lambda}\ln\left(\frac{r-k}{c}\right).
    \intertext{Solving $\hat{\pi}_s(\rho)=\reservationprofit$ for $\rho$ gives:}
    \hat{\rho} &= \frac{1}{\lambda}\left(r-k-c-\frac{c}{\bb\lambda+1}\ln\left(\frac{r-k}{c}\right)-\reservationprofit\frac{\lambda}{b\lambda+1}\right).
    \intertext{Note that $\hat{\rho}>0$ whenever $r-k-c>\frac{c}{\bb\lambda+1}\ln\left(\frac{r-k}{c}\right)+\reservationprofit\frac{\lambda}{b\lambda+1}$, or equivalently $\reservationprofit < (r-k-c)\frac{\bb\lambda+1}{\lambda}-\frac{c}{\lambda}\ln\left(\frac{r-k}{c}\right)=\Pi^*$. To satisfy $r=w+\rho\lambda$, the OEM sets}
    \hat{w}&=k+c+\frac{c}{\bb\lambda+1}\ln\left(\frac{r-k}{c}\right)+\reservationprofit\frac{\lambda}{\bb\lambda+1}.
    \intertext{By construction $\hat{\pi}_s(\hat{x}(\hat{w},\hat{\rho}),\hat{w},\hat{\rho})=\reservationprofit$ and it follows that}
    \hat{\pi}_m(\hat{w},\hat{\rho}) &= \left(r-k-c-\frac{c}{\bb\lambda+1}\ln\left(\frac{r-k}{c}\right)\right)\left(\bb+\frac{1}{\lambda}\right)-\reservationprofit\\
    &= \Pi^*-\reservationprofit.
\end{align*}
\end{proof}

\paragraph{Proof of Lemma \ref{lem: extension - exp sales}}
\begin{proof}
Since,
\begin{align*}
    \mathbb{E}\left[ \min\{D,x\} \right] &= \mathbb{E}\left[ D | D<x\right]\mathbb{P}(D<x) + x\mathbb{P}(D>x)
\end{align*}
we need to determine $\mathbb{E}\left[ D | D<x\right]$.
\begin{align*}
    \mathbb{E}\left[ D | D<x\right]= \frac{\int_0^x y \frac{\lambda^\erlangshape y^{\erlangshape-1}}{(\erlangshape-1)!}e^{-\lambda y}dy}{\mathbb{P}(D<x)}
\end{align*}
The lower incomplete gamma function is defined as:
\begin{align*}
    \gamma(s,x)&=\int_0^x t^{s-1}e^{-t} dt
    \intertext{such that}
    \gamma(\erlangshape+1,\lambda x)&=\int_0^x t^\erlangshape e^{-t} dt
    \intertext{and}
    \frac{1}{\lambda}\gamma(\erlangshape+1,\lambda x)&=\int_0^x (\lambda y)^\erlangshape e^{-\lambda y} dy.
\end{align*}
Therefore, 
\begin{align*}
    \int_0^x y \frac{\lambda^\erlangshape y^{\erlangshape-1}}{(\erlangshape-1)!}e^{-\lambda y}dy &= \frac{1}{\lambda(\erlangshape-1)!}\gamma(\erlangshape+1,\lambda x)
    \intertext{and consequently,}
    \mathbb{E}\left[ \min\{D,x\} \right]&=\frac{\gamma(k+1,\lambda x)}{\lambda(\erlangshape-1)!} + x \mathbb{P}(D>x).
\end{align*}
\end{proof}

\end{appendices}